\documentclass[12pt,a4paper]{article}
\usepackage{amsmath,amsthm,amssymb,latexsym,epic}
\usepackage{bezier}

\newtheorem{theorem}{Theorem}
\newtheorem{lemma}{Lemma}
\newtheorem{corollary}{Corollary}
\newtheorem{proposition}{Proposition}

\usepackage[all]{xy}
\usepackage{enumerate}

\newcommand{\IP}{\mathcal{IP}}
\newcommand{\PIP}{\mathcal{PIP}}

\newcommand{\Sym}{\mathcal{S}}

\renewcommand{\S}{\mathcal{S}}

\newcommand{\C}{\mathfrak{C}}

\renewcommand{\phi}{\varphi}

\begin{document}

\title{A presentation for the partial dual inverse symmetric monoid}
\author{Ganna Kudryavtseva and Victor Maltcev}
\date{}
\maketitle

\begin{abstract}
We give a monoid presentation in terms of generators and defining
relations for the partial analogue of the finite dual inverse
symmetric monoid.
\end{abstract}

AMS 2000 Mathematical Subject Classification 20M05, 68R15, 20M20.

\section{Introduction}\label{s1}

The partial dual inverse symmetric monoid on a set $X$, denoted by
$\PIP_X$, is a partial analogue of the dual inverse symmetric
monoid $\IP_X$, see~\cite{M} and~\cite{FL}. This monoid is a
natural generalization of the full inverse symmetric monoid and
has a number of interesting properties, which were studied
in~\cite{KMal}.

The aim of the present paper is to obtain a presentation for the
monoid $\PIP_X$ with $X$ finite in terms of generators and
defining relations. We would like to mention that during the
recent period there appeared a number of papers where
presentations for some important transformation semigroups and
their generalizations (e.g., so called Brauer-type semigroups)
have been found, see~\cite{Fer}, \cite{F}, \cite{E}, \cite{KMaz},
\cite{MM}. In view of this, our research looks like a natural
continuation of the previous efforts.

It is interesting that by the moment (so far as to our knowledge)
a presentation for the finite dual inverse symmetric monoid
$\IP_X$ is not found. Some possible approaches towards finding
such a presentation and the arising difficulties are discussed
in~\cite{EEF}. In view of this, our result looks somehow
unexpected, as we solve the problem for a bigger and more
complicated monoid. The authors have a hope that the ideas and
technique suggested in the present paper could be utilized, in
particular, for finding a presentation for the monoid $\IP_X$.

The paper is organized as follows. In Section~\ref{sec:def-PIP} we
recall the definition of the monoid $\PIP_X$ (and of the monoid
$\IP_X$). In Section~\ref{s3} we define an abstract monoid $S$ by
generators and defining relations and establish some other
relations which are the consequences of the defining ones.
Further, in Section~\ref{s4} we continue investigating the monoid
$S$ and develop some rewriting technique for the elements of $S$
presented as words over its generators. Using this technique we
manage to show that every element of $S$ can be presented as a
certain "canonical word". Finally, in Section~\ref{s5} we turn
back to the monoid $\PIP_X$, $X$ finite, construct a natural
epimorphism from $S$ onto $\PIP_X$ and show that this epimorphism
is in fact an isomorphism. For this, we prove that the
presentation of an element of $\PIP_X$ as an image of some
canonical word is unique, and thus the cardinality of $S$ does not
exceed the cardinality of $\PIP_X$.

\section{Definition of the monoid $\PIP_n$}\label{sec:def-PIP}

Let $X$ be a set. Consider a set $X'=\{x'\}_{x\in X}$ disjoint
with $X$ and a bijection $':X\to X'$ sending $x\in X$ to $x'\in
X'$. Denote the inverse bijection by the same symbol, that is
$(x')'=x$ for all $x\in X\cup X'$.

We shall say that a subset $A$ of $X\cup X'$ is a
\begin{itemize}
\item
\emph{line} provided that $A\cap X\ne\varnothing$ and $A\cap
X'\ne\varnothing$;
\item
\emph{point} provided that $\mid A\mid=1$.
\end{itemize}

Let $\IP_X$  be the set of all decompositions of $X\cup X'$ into
lines, and $\PIP_X$  the set of all decompositions of $X\cup X'$
into lines and points. Obviously, $\IP_X\subset \PIP_X$.

In the case when $X=\{1,\ldots,n\}$ we shall denote $\IP_X$ by
$\IP_n$ and $\PIP_X$ by $\PIP_n$.

Let $a\in\PIP_X$ and $x,y\in X\cup X'$. Set $x\equiv_{a}y$
provided that $x$ and $y$ are of the same block of $a$. The map
$a\mapsto \equiv_a$ is a bijection between the elements of
$\PIP_X$ and the equivalence relations on $X\cup X'$ whose classes
are either points or lines. Under this bijection the set $\IP_X$
maps onto the set of those equivalence relations on $X\cup X'$
whose classes are lines.

To define a multiplication on the set $\IP_X$ we consider any
$a,b\in\IP_X$ and define a new equivalence relation, $\equiv$, on
$X\cup X'$ as follows:
\begin{itemize}
\item
for $x,y\in X$ we have $x\equiv y$ if and only if $x\equiv_{a} y$
or there is a sequence, $c_1,\ldots,c_{2s}$, $s\geq 1$, of
elements of $X$, such that $x\equiv_{a} c'_1$, $c_1\equiv_{b}
c_2$, $c'_2\equiv_{a} c'_3,\ldots,$ $c_{2s-1}\equiv_{b} c_{2s}$,
and $c'_{2s}\equiv_{a} y$;
\item
for $x,y\in X$ we have $x'\equiv y'$ if and only if $x'\equiv_{b}
y'$ or there is a sequence, $c_1,\ldots,c_{2s}$, $s\geq 1$, of
elements of $X$, such that  $x'\equiv_{b} c_1$, $c'_1\equiv_{a}
c'_2$, $c_2\equiv_{b} c_3,\ldots,$ $c'_{2s-1}\equiv_{a} c'_{2s}$,
and $c_{2s}\equiv_{b} y'$;
\item
for $x,y\in X$ we have $x\equiv y'$ if and only if $y'\equiv x$ if
and only if there is a sequence, $c_1,\ldots$, $c_{2s-1}$, $s\geq
1$, of elements of $X$, such that $x\equiv_{a} c'_1$,
$c_1\equiv_{b} c_2$, $c'_2\equiv_{a} c'_3,\ldots,$
$c'_{2s-2}\equiv_{a} c'_{2s-1}$, and $c_{2s-1}\equiv_{b} y'$.
\end{itemize}

Since every class of $\equiv$ is a line this definition is
correct. We set the decomposition of $X\cup X'$ into
$\equiv$-classes to be the product $a\cdot b$ of $a$ and $b$ in
$\IP_X$. With respect to this multiplication
$\bigl(\IP_X,\cdot\bigr)$ is a semigroup. It was called the
\emph{inverse partition semigroup} on the set $X$ in \cite{M}
and~\cite{M1}, the {\em monoid of block bijections} in~\cite{EEF},
and the {\em dual inverse symmetric monoid} on the set $X$
in~\cite{FL} and a number of subsequent papers. In this paper we
stick to the latter term.

Let $x\not\in X$ be an arbitrary element. Set $Y=X\cup \{x\}$ and
denote by $\widetilde{\IP}_{Y}$ the subset of $\IP_{Y}$ consisting
of those decomposition of $Y\cup Y'$ into subsets which consist
entirely of lines and both $x$ and $x'$ belong to the same line.
The set $\widetilde{\IP}_{Y}$ is closed with respect to the
operation $\cdot$ and is therefore a  subsemigroup of $\IP_{Y}$.

Take $a\in\PIP_X$ and denote by $\varphi(a)$ the element of
$\widetilde{\IP}_{Y}$, consisting of all lines of $a$ and of one
additional block, whose elements are $x, x'$ and all points of
$a$. It was noticed in~\cite{KMal} (and is easy to see) that the
map $\varphi$ is a bijection from the set $\PIP_X$ onto the set
$\widetilde{\IP}_{Y}$. Now we are prepared to define the
(associative) multiplication on $\PIP_X$. We set (slightly abusing
the notation)
\begin{equation*}
a\cdot b=\varphi^{-1}\bigl(\varphi(a)\cdot\varphi(b)\bigr).
\end{equation*}

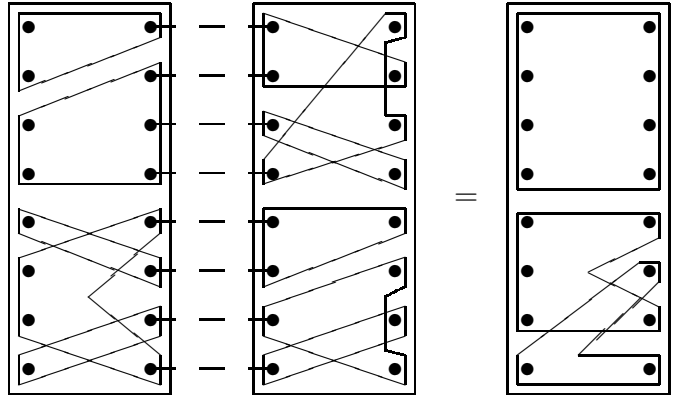
\begin{figure}
\special{em:linewidth 0.4pt} \unitlength 0.65mm
\linethickness{0.4pt}
\begin{picture}(150.00,80.00)(-20.00,0.00)
\put(20.00,00.00){\makebox(0,0)[cc]{$\bullet$}}
\put(20.00,10.00){\makebox(0,0)[cc]{$\bullet$}}
\put(20.00,20.00){\makebox(0,0)[cc]{$\bullet$}}
\put(20.00,30.00){\makebox(0,0)[cc]{$\bullet$}}
\put(20.00,40.00){\makebox(0,0)[cc]{$\bullet$}}
\put(20.00,50.00){\makebox(0,0)[cc]{$\bullet$}}
\put(20.00,60.00){\makebox(0,0)[cc]{$\bullet$}}
\put(20.00,70.00){\makebox(0,0)[cc]{$\bullet$}}
\put(45.00,10.00){\makebox(0,0)[cc]{$\bullet$}}
\put(45.00,20.00){\makebox(0,0)[cc]{$\bullet$}}
\put(45.00,30.00){\makebox(0,0)[cc]{$\bullet$}}
\put(45.00,40.00){\makebox(0,0)[cc]{$\bullet$}}
\put(45.00,50.00){\makebox(0,0)[cc]{$\bullet$}}
\put(45.00,60.00){\makebox(0,0)[cc]{$\bullet$}}
\put(45.00,70.00){\makebox(0,0)[cc]{$\bullet$}}
\put(45.00,00.00){\makebox(0,0)[cc]{$\bullet$}}

\drawline(16.00,-05.00)(16.00,75.00)
\drawline(16.00,75.00)(49.00,75.00)
\drawline(49.00,75.00)(49.00,-05.00)
\drawline(49.00,-05.00)(16.00,-05.00)

\drawline(18.00,-03.00)(18.00,03.00)
\drawline(18.00,03.00)(47.00,13.00)
\drawline(47.00,13.00)(47.00,07.00)
\drawline(47.00,07.00)(18.00,-03.00)
\drawline(18.00,07.00)(18.00,23.00)
\drawline(18.00,23.00)(47.00,33.00)
\drawline(47.00,33.00)(47.00,28.00)
\drawline(47.00,28.00)(32.25,15.00)
\drawline(32.25,15.00)(47.00,03.00)
\drawline(47.00,03.00)(47.00,-03.00)
\drawline(47.00,-03.00)(18.00,07.00)
\drawline(18.00,28.00)(18.00,33.00)
\drawline(18.00,33.00)(47.00,23.00)
\drawline(47.00,23.00)(47.00,17.00)
\drawline(47.00,17.00)(18.00,28.00)
\drawline(18.00,38.00)(18.00,52.00)
\drawline(18.00,52.00)(47.00,63.00)
\drawline(47.00,63.00)(47.00,38.00)
\drawline(47.00,38.00)(18.00,38.00)
\drawline(18.00,57.00)(18.00,73.00)
\drawline(18.00,73.00)(47.00,73.00)
\drawline(47.00,73.00)(47.00,68.00)
\drawline(47.00,68.00)(18.00,57.00)
\put(70.00,00.00){\makebox(0,0)[cc]{$\bullet$}}
\put(70.00,10.00){\makebox(0,0)[cc]{$\bullet$}}
\put(70.00,20.00){\makebox(0,0)[cc]{$\bullet$}}
\put(70.00,30.00){\makebox(0,0)[cc]{$\bullet$}}
\put(70.00,40.00){\makebox(0,0)[cc]{$\bullet$}}
\put(70.00,50.00){\makebox(0,0)[cc]{$\bullet$}}
\put(70.00,60.00){\makebox(0,0)[cc]{$\bullet$}}
\put(70.00,70.00){\makebox(0,0)[cc]{$\bullet$}}
\put(95.00,10.00){\makebox(0,0)[cc]{$\bullet$}}
\put(95.00,20.00){\makebox(0,0)[cc]{$\bullet$}}
\put(95.00,30.00){\makebox(0,0)[cc]{$\bullet$}}
\put(95.00,40.00){\makebox(0,0)[cc]{$\bullet$}}
\put(95.00,50.00){\makebox(0,0)[cc]{$\bullet$}}
\put(95.00,60.00){\makebox(0,0)[cc]{$\bullet$}}
\put(95.00,70.00){\makebox(0,0)[cc]{$\bullet$}}
\put(95.00,00.00){\makebox(0,0)[cc]{$\bullet$}}

\drawline(66.00,-05.00)(66.00,75.00)
\drawline(66.00,75.00)(99.00,75.00)
\drawline(99.00,75.00)(99.00,-05.00)
\drawline(99.00,-05.00)(66.00,-05.00)

\drawline(68.00,-03.00)(68.00,03.00)
\drawline(68.00,03.00)(97.00,13.00)
\drawline(97.00,13.00)(97.00,07.00)
\drawline(97.00,07.00)(68.00,-03.00)
\drawline(68.00,07.00)(68.00,13.00)
\drawline(68.00,13.00)(97.00,23.00)
\drawline(97.00,23.00)(97.00,17.00)
\drawline(97.00,17.00)(93.00,15.00)
\drawline(93.00,15.00)(93.00,04.00)
\drawline(93.00,04.00)(97.00,03.00)
\drawline(97.00,03.00)(97.00,-03.00)
\drawline(97.00,-03.00)(68.00,07.00)
\drawline(68.00,17.00)(68.00,33.00)
\drawline(68.00,33.00)(97.00,33.00)
\drawline(97.00,33.00)(97.00,28.00)
\drawline(97.00,28.00)(68.00,17.00)
\drawline(68.00,48.00)(68.00,53.00)
\drawline(68.00,53.00)(97.00,43.00)
\drawline(97.00,43.00)(97.00,37.00)
\drawline(97.00,37.00)(68.00,48.00)
\drawline(68.00,58.00)(68.00,73.00)
\drawline(68.00,73.00)(97.00,63.00)
\drawline(97.00,63.00)(97.00,58.00)
\drawline(97.00,58.00)(68.00,58.00)
\drawline(68.00,38.00)(68.00,43.00)
\drawline(68.00,43.00)(93.00,73.00)
\drawline(93.00,73.00)(97.00,73.00)
\drawline(97.00,73.00)(97.00,68.00)
\drawline(97.00,68.00)(93.00,67.00)
\drawline(93.00,67.00)(93.00,52.00)
\drawline(93.00,52.00)(97.00,52.00)
\drawline(97.00,52.00)(97.00,47.00)
\drawline(97.00,47.00)(68.00,38.00)

\drawline(45.00,00.30)(50.00,00.30)
\drawline(55.00,00.30)(60.00,00.30)
\drawline(65.00,00.30)(70.00,00.30)
\drawline(45.00,10.30)(50.00,10.30)
\drawline(55.00,10.30)(60.00,10.30)
\drawline(65.00,10.30)(70.00,10.30)
\drawline(45.00,20.30)(50.00,20.30)
\drawline(55.00,20.30)(60.00,20.30)
\drawline(65.00,20.30)(70.00,20.30)
\drawline(45.00,30.30)(50.00,30.30)
\drawline(55.00,30.30)(60.00,30.30)
\drawline(65.00,30.30)(70.00,30.30)
\drawline(45.00,40.30)(50.00,40.30)
\drawline(55.00,40.30)(60.00,40.30)
\drawline(65.00,40.30)(70.00,40.30)
\drawline(45.00,50.30)(50.00,50.30)
\drawline(55.00,50.30)(60.00,50.30)
\drawline(65.00,50.30)(70.00,50.30)
\drawline(45.00,60.30)(50.00,60.30)
\drawline(55.00,60.30)(60.00,60.30)
\drawline(65.00,60.30)(70.00,60.30)
\drawline(45.00,70.30)(50.00,70.30)
\drawline(55.00,70.30)(60.00,70.30)
\drawline(65.00,70.30)(70.00,70.30)
\drawline(120.00,37.00)(120.00,73.00)
\drawline(120.00,73.00)(149.00,73.00)
\drawline(149.00,73.00)(149.00,37.00)
\drawline(149.00,37.00)(120.00,37.00)
\drawline(120.00,08.00)(120.00,32.00)
\drawline(120.00,32.00)(149.00,32.00)
\drawline(149.00,32.00)(149.00,27.00)
\drawline(149.00,27.00)(134.50,20.00)
\drawline(134.50,20.00)(149.00,13.00)
\drawline(149.00,13.00)(149.00,08.00)
\drawline(149.00,08.00)(120.00,08.00)
\drawline(120.00,-03.00)(120.00,03.00)
\drawline(120.00,03.00)(145.00,22.00)
\drawline(145.00,22.00)(149.00,22.00)
\drawline(149.00,22.00)(149.00,18.00)
\drawline(149.00,18.00)(132.50,03.00)
\drawline(132.50,03.00)(149.00,03.00)
\drawline(149.00,03.00)(149.00,-03.00)
\drawline(149.00,-03.00)(120.00,-03.00)
\put(122.00,00.00){\makebox(0,0)[cc]{$\bullet$}}
\put(122.00,10.00){\makebox(0,0)[cc]{$\bullet$}}
\put(122.00,20.00){\makebox(0,0)[cc]{$\bullet$}}
\put(122.00,30.00){\makebox(0,0)[cc]{$\bullet$}}
\put(122.00,40.00){\makebox(0,0)[cc]{$\bullet$}}
\put(122.00,50.00){\makebox(0,0)[cc]{$\bullet$}}
\put(122.00,60.00){\makebox(0,0)[cc]{$\bullet$}}
\put(122.00,70.00){\makebox(0,0)[cc]{$\bullet$}}
\put(147.00,10.00){\makebox(0,0)[cc]{$\bullet$}}
\put(147.00,20.00){\makebox(0,0)[cc]{$\bullet$}}
\put(147.00,30.00){\makebox(0,0)[cc]{$\bullet$}}
\put(147.00,40.00){\makebox(0,0)[cc]{$\bullet$}}
\put(147.00,50.00){\makebox(0,0)[cc]{$\bullet$}}
\put(147.00,60.00){\makebox(0,0)[cc]{$\bullet$}}
\put(147.00,70.00){\makebox(0,0)[cc]{$\bullet$}}
\put(147.00,00.00){\makebox(0,0)[cc]{$\bullet$}}
\put(109.50,35.00){\makebox(0,0)[cc]{$=$}}

\drawline(118.00,-05.00)(118.00,75.00)
\drawline(118.00,75.00)(151.00,75.00)
\drawline(151.00,75.00)(151.00,-05.00)
\drawline(151.00,-05.00)(118.00,-05.00)
\end{picture}
\caption{Elements of $\IP_8$ and their
multiplication.}\label{fig:ip}
\end{figure}

\begin{figure}
\special{em:linewidth 0.4pt} \unitlength 0.65mm
\linethickness{0.4pt}
\begin{picture}(150.00,90.00)(-20.00, 0.00)
\put(20.00,00.00){\makebox(0,0)[cc]{$\bullet$}}
\put(20.00,10.00){\makebox(0,0)[cc]{$\bullet$}}
\put(20.00,20.00){\makebox(0,0)[cc]{$\bullet$}}
\put(20.00,30.00){\makebox(0,0)[cc]{$\bullet$}}
\put(20.00,40.00){\makebox(0,0)[cc]{$\bullet$}}
\put(20.00,50.00){\makebox(0,0)[cc]{$\bullet$}}
\put(20.00,60.00){\makebox(0,0)[cc]{$\bullet$}}
\put(20.00,70.00){\makebox(0,0)[cc]{$\bullet$}}
\put(45.00,10.00){\makebox(0,0)[cc]{$\bullet$}}
\put(45.00,20.00){\makebox(0,0)[cc]{$\bullet$}}
\put(45.00,30.00){\makebox(0,0)[cc]{$\bullet$}}
\put(45.00,40.00){\makebox(0,0)[cc]{$\bullet$}}
\put(45.00,50.00){\makebox(0,0)[cc]{$\bullet$}}
\put(45.00,60.00){\makebox(0,0)[cc]{$\bullet$}}
\put(45.00,70.00){\makebox(0,0)[cc]{$\bullet$}}
\put(45.00,00.00){\makebox(0,0)[cc]{$\bullet$}}
\drawline(16.00,-05.00)(16.00,75.00)
\drawline(16.00,75.00)(49.00,75.00)
\drawline(49.00,75.00)(49.00,-05.00)
\drawline(49.00,-05.00)(16.00,-05.00)
\drawline(66.00,-05.00)(66.00,75.00)
\drawline(66.00,75.00)(99.00,75.00)
\drawline(99.00,75.00)(99.00,-05.00)
\drawline(99.00,-05.00)(66.00,-05.00)
\drawline(118.00,-05.00)(118.00,75.00)
\drawline(118.00,75.00)(151.00,75.00)
\drawline(151.00,75.00)(151.00,-05.00)
\drawline(151.00,-05.00)(118.00,-05.00)

\drawline(18.00,-03.00)(18.00,13.00)
\drawline(18.00,13.00)(22.00,13.00)
\drawline(22.00,13.00)(43.00,23.00)
\drawline(43.00,23.00)(47.00,23.00)
\drawline(47.00,23.00)(47.00,-03.00)
\drawline(47.00,-03.00)(18.00,-03.00)

\drawline(18.00,18.00)(18.00,32.00)
\drawline(18.00,32.00)(47.00,32.00)
\drawline(47.00,32.00)(47.00,27.00)
\drawline(47.00,27.00)(22.00,18.00)
\drawline(22.00,18.00)(18.00,18.00)

\drawline(18.00,38.00)(18.00,42.00)
\drawline(18.00,42.00)(43.00,52.00)
\drawline(43.00,52.00)(47.00,52.00)
\drawline(47.00,52.00)(47.00,48.00)
\drawline(47.00,48.00)(22.00,38.00)
\drawline(22.00,38.00)(18.00,38.00)

\drawline(18.00,48.00)(18.00,62.00)
\drawline(18.00,62.00)(43.00,72.00)
\drawline(43.00,72.00)(47.00,72.00)
\drawline(47.00,72.00)(47.00,58.00)
\drawline(47.00,58.00)(22.00,48.00)
\drawline(22.00,48.00)(18.00,48.00)

\drawline(18.00,68.00)(18.00,72.00)
\drawline(18.00,72.00)(22.00,72.00)
\drawline(22.00,72.00)(22.00,68.00)
\drawline(22.00,68.00)(18.00,68.00)

\drawline(43.00,38.00)(43.00,42.00)
\drawline(43.00,42.00)(47.00,42.00)
\drawline(47.00,42.00)(47.00,38.00)
\drawline(47.00,38.00)(43.00,38.00)

\drawline(68.00,-03.00)(68.00,03.00)
\drawline(68.00,03.00)(91.00,22.00)
\drawline(91.00,22.00)(97.00,22.00)
\drawline(97.00,22.00)(97.00,18.00)
\drawline(97.00,18.00)(91.00,18.00)
\drawline(91.00,18.00)(91.00,03.00)
\drawline(91.00,03.00)(97.00,03.00)
\drawline(97.00,03.00)(97.00,-03.00)
\drawline(97.00,-03.00)(68.00,-03.00)

\drawline(68.00,08.00)(68.00,23.00)
\drawline(68.00,23.00)(93.00,32.00)
\drawline(93.00,32.00)(97.00,32.00)
\drawline(97.00,32.00)(97.00,28.00)
\drawline(97.00,28.00)(93.00,28.00)
\drawline(93.00,28.00)(72.00,08.00)
\drawline(72.00,08.00)(68.00,08.00)

\drawline(68.00,28.00)(68.00,32.00)
\drawline(68.00,32.00)(74.00,32.00)
\drawline(74.00,32.00)(74.00,48.00)
\drawline(74.00,48.00)(68.00,48.00)
\drawline(68.00,48.00)(68.00,52.00)
\drawline(68.00,52.00)(97.00,52.00)
\drawline(97.00,52.00)(97.00,38.00)
\drawline(97.00,38.00)(93.00,38.00)
\drawline(93.00,38.00)(72.00,28.00)
\drawline(72.00,28.00)(68.00,28.00)

\drawline(68.00,68.00)(68.00,72.00)
\drawline(68.00,72.00)(97.00,72.00)
\drawline(97.00,72.00)(97.00,58.00)
\drawline(97.00,58.00)(93.00,58.00)
\drawline(93.00,58.00)(68.00,68.00)

\drawline(68.00,38.00)(68.00,42.00)
\drawline(68.00,42.00)(72.00,42.00)
\drawline(72.00,42.00)(72.00,38.00)
\drawline(72.00,38.00)(68.00,38.00)

\drawline(68.00,58.00)(68.00,62.00)
\drawline(68.00,62.00)(72.00,62.00)
\drawline(72.00,62.00)(72.00,58.00)
\drawline(72.00,58.00)(68.00,58.00)

\drawline(120.00,-03.00)(120.00,12.00)
\drawline(120.00,12.00)(124.00,12.00)
\drawline(124.00,12.00)(145.00,32.00)
\drawline(145.00,32.00)(149.00,32.00)
\drawline(149.00,32.00)(149.00,18.00)
\drawline(149.00,18.00)(143.00,18.00)
\drawline(143.00,18.00)(143.00,02.00)
\drawline(143.00,02.00)(149.00,02.00)
\drawline(149.00,02.00)(149.00,-03.00)
\drawline(149.00,-03.00)(120.00,-03.00)

\drawline(120.00,18.00)(120.00,42.00)
\drawline(120.00,42.00)(145.00,52.00)
\drawline(145.00,52.00)(149.00,52.00)
\drawline(149.00,52.00)(149.00,38.00)
\drawline(149.00,38.00)(145.00,38.00)
\drawline(145.00,38.00)(124.00,18.00)
\drawline(124.00,18.00)(120.00,18.00)

\drawline(120.00,48.00)(120.00,52.00)
\drawline(120.00,52.00)(124.00,52.00)
\drawline(124.00,52.00)(124.00,48.00)
\drawline(124.00,48.00)(120.00,48.00)

\drawline(120.00,58.00)(120.00,62.00)
\drawline(120.00,62.00)(124.00,62.00)
\drawline(124.00,62.00)(124.00,58.00)
\drawline(124.00,58.00)(120.00,58.00)

\drawline(120.00,68.00)(120.00,72.00)
\drawline(120.00,72.00)(124.00,72.00)
\drawline(124.00,72.00)(124.00,68.00)
\drawline(124.00,68.00)(120.00,68.00)

\drawline(145.00,08.00)(145.00,12.00)
\drawline(145.00,12.00)(149.00,12.00)
\drawline(149.00,12.00)(149.00,08.00)
\drawline(149.00,08.00)(145.00,08.00)

\drawline(145.00,58.00)(145.00,62.00)
\drawline(145.00,62.00)(149.00,62.00)
\drawline(149.00,62.00)(149.00,58.00)
\drawline(149.00,58.00)(145.00,58.00)

\drawline(145.00,68.00)(145.00,72.00)
\drawline(145.00,72.00)(149.00,72.00)
\drawline(149.00,72.00)(149.00,68.00)
\drawline(149.00,68.00)(145.00,68.00)

\drawline(93.00,08.00)(93.00,12.00)
\drawline(93.00,12.00)(97.00,12.00)
\drawline(97.00,12.00)(97.00,08.00)
\drawline(97.00,08.00)(93.00,08.00)

\put(70.00,00.00){\makebox(0,0)[cc]{$\bullet$}}
\put(70.00,10.00){\makebox(0,0)[cc]{$\bullet$}}
\put(70.00,20.00){\makebox(0,0)[cc]{$\bullet$}}
\put(70.00,30.00){\makebox(0,0)[cc]{$\bullet$}}
\put(70.00,40.00){\makebox(0,0)[cc]{$\bullet$}}
\put(70.00,50.00){\makebox(0,0)[cc]{$\bullet$}}
\put(70.00,60.00){\makebox(0,0)[cc]{$\bullet$}}
\put(70.00,70.00){\makebox(0,0)[cc]{$\bullet$}}
\put(95.00,10.00){\makebox(0,0)[cc]{$\bullet$}}
\put(95.00,20.00){\makebox(0,0)[cc]{$\bullet$}}
\put(95.00,30.00){\makebox(0,0)[cc]{$\bullet$}}
\put(95.00,40.00){\makebox(0,0)[cc]{$\bullet$}}
\put(95.00,50.00){\makebox(0,0)[cc]{$\bullet$}}
\put(95.00,60.00){\makebox(0,0)[cc]{$\bullet$}}
\put(95.00,70.00){\makebox(0,0)[cc]{$\bullet$}}
\put(95.00,00.00){\makebox(0,0)[cc]{$\bullet$}}
\put(122.00,00.00){\makebox(0,0)[cc]{$\bullet$}}
\put(122.00,10.00){\makebox(0,0)[cc]{$\bullet$}}
\put(122.00,20.00){\makebox(0,0)[cc]{$\bullet$}}
\put(122.00,30.00){\makebox(0,0)[cc]{$\bullet$}}
\put(122.00,40.00){\makebox(0,0)[cc]{$\bullet$}}
\put(122.00,50.00){\makebox(0,0)[cc]{$\bullet$}}
\put(122.00,60.00){\makebox(0,0)[cc]{$\bullet$}}
\put(122.00,70.00){\makebox(0,0)[cc]{$\bullet$}}
\put(147.00,10.00){\makebox(0,0)[cc]{$\bullet$}}
\put(147.00,20.00){\makebox(0,0)[cc]{$\bullet$}}
\put(147.00,30.00){\makebox(0,0)[cc]{$\bullet$}}
\put(147.00,40.00){\makebox(0,0)[cc]{$\bullet$}}
\put(147.00,50.00){\makebox(0,0)[cc]{$\bullet$}}
\put(147.00,60.00){\makebox(0,0)[cc]{$\bullet$}}
\put(147.00,70.00){\makebox(0,0)[cc]{$\bullet$}}
\put(147.00,00.00){\makebox(0,0)[cc]{$\bullet$}}
\drawline(45.00,00.30)(50.00,00.30)
\drawline(55.00,00.30)(60.00,00.30)
\drawline(65.00,00.30)(70.00,00.30)
\drawline(45.00,10.30)(50.00,10.30)
\drawline(55.00,10.30)(60.00,10.30)
\drawline(65.00,10.30)(70.00,10.30)
\drawline(45.00,20.30)(50.00,20.30)
\drawline(55.00,20.30)(60.00,20.30)
\drawline(65.00,20.30)(70.00,20.30)
\drawline(45.00,30.30)(50.00,30.30)
\drawline(55.00,30.30)(60.00,30.30)
\drawline(65.00,30.30)(70.00,30.30)
\drawline(45.00,40.30)(50.00,40.30)
\drawline(55.00,40.30)(60.00,40.30)
\drawline(65.00,40.30)(70.00,40.30)
\drawline(45.00,50.30)(50.00,50.30)
\drawline(55.00,50.30)(60.00,50.30)
\drawline(65.00,50.30)(70.00,50.30)
\drawline(45.00,60.30)(50.00,60.30)
\drawline(55.00,60.30)(60.00,60.30)
\drawline(65.00,60.30)(70.00,60.30)
\drawline(45.00,70.30)(50.00,70.30)
\drawline(55.00,70.30)(60.00,70.30)
\drawline(65.00,70.30)(70.00,70.30)
\put(109.50,35.00){\makebox(0,0)[cc]{$=$}}
\end{picture}
\caption{Elements of $\PIP_8$ and their
multiplication.}\label{fig:pip}

\end{figure}
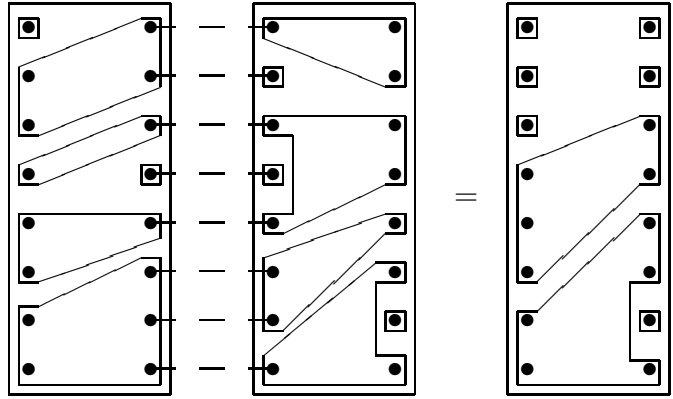

The above defined multiplication in the monoid $\PIP_X$ has a
natural realization as a "superposition of diagrams".  We
interpret the elements of $\PIP_X$ as diagrams with vertices on
the left hand side indexed by $X$ and vertices on the right hand
side indexed by $X'$. To multiply two such diagrams $\alpha$ and
$\beta$ one has to place $\beta$ to the right of $\alpha$ such
that the corresponding right vertices of $\alpha$ and left
vertices of $\beta$ are identified, which uniquely determines the
diagram of the product decomposition $\alpha\beta$. This is
illustrated on Figures~\ref{fig:ip} and~\ref{fig:pip}.

The semigroup $\bigl(\PIP_X,\cdot\bigr)$ is a "partial analogue"
of the semigroup $\IP_X$. In particular, it contains the semigroup
$\bigl(\IP_X,\cdot\bigr)$ as a subsemigroup. The structure of the
semigroup $\bigl(\PIP_X,\cdot\bigr)$ was investigated in
\cite{KMal}, where it was called the {\em partial inverse
partition semigroup}. We would like, developing the terminology
stemming from~\cite{FL}, to propose a more apt, from our point of
view, term for the monoid $\PIP_X$, the {\em partial dual inverse
symmetric monoid}.

\section{The abstract monoid $S$ and its relations}\label{s3}

Let $n\geq 3$. Consider the monoid $S$ with the identity element
$e$ generated by $\sigma_1, \dots, \sigma_{n-1}$;
$\lambda_1,\dots, \lambda_{n-1}$; $\rho_1, \dots, \rho_{n-1}$;
$e_1,\dots, e_n$ subject to the following relations:
\begin{gather}
\sigma_i^2=e, 1\leq i\leq n-1, \label{c1}\\
\sigma_i\sigma_j=\sigma_j\sigma_i, |i-j|>1, \label{c2}\\
\sigma_i\sigma_j\sigma_i=\sigma_j\sigma_i\sigma_j, |i-j|=1, \label{c3} \\
\lambda_i\lambda_j=\lambda_j\lambda_i, \, \rho_i\rho_j=\rho_j\rho_i, \,\,
\rho_i\lambda_j=\lambda_j\rho_i, \, |i-j|>1, \label{lr1}\\
\lambda_i\rho_{i+1}=\sigma_{i+1}\rho_i\lambda_{i}e_{i+2}=\sigma_i\rho_{i+1}\sigma_i\lambda_i,
\,\,
\lambda_{i+1}\rho_i=e_{i+2}\rho_{i}\lambda_i\sigma_{i+1}=\rho_i\sigma_i\lambda_{i+1}\sigma_i,\label{neighb-lambda-rho}\\
\rho_i\lambda_{i+1}=\rho_ie_{i+2}, \,\,
\rho_{i+1}\lambda_i=e_{i+2}\lambda_i, \label{neighb-rho-lambda}\\
\lambda_i\rho_i\lambda_i=\lambda_i, \rho_i\lambda_i\rho_i=\rho_i,
\label{inverse}\\
 \sigma_i\lambda_j\sigma_i=\sigma_j\lambda_i\sigma_j,\,\,
\sigma_i\rho_j\sigma_i=\sigma_j\rho_i\sigma_j,|i-j|=1,
\label{slr1}\\
\lambda_i\sigma_i=\lambda_i, \sigma_i\rho_i=\rho_i, 1\leq i\leq
n-1, \label{slr2}\\
\sigma_i\lambda_j=\lambda_j\sigma_i, \,
\sigma_i\rho_j=\rho_j\sigma_i,\, |j-i|>1, \label{slr3}\\
\lambda_i\lambda_{i+1}=\lambda_{i}\lambda_{i+1}\sigma_i=\sigma_{i+1}\lambda_{i}\lambda_{i+1},\,\
\rho_{i+1}\rho_{i}=\sigma_{i}\rho_{i+1}\rho_{i}=\rho_{i+1}\rho_{i}\sigma_{i+1},\label{ghost}\\
\sigma_{i+1}\lambda_{i+1}\lambda_i=\lambda_{i+1}\lambda_i=\lambda_ie_{i+2},\,\,
\rho_i\rho_{i+1}\sigma_{i+1}=\rho_i\rho_{i+1}= e_{i+2}\rho_i,\label{ghost1}\\
 e_i=\lambda_{i-1}\rho_{i-1}, i\geq 2, \,\,
e_1=\sigma_1e_2\sigma_1,
\label{df-e}\\
e_i^2=e_i,
e_ie_{i+1}=e_{i+1}e_i=\lambda_i^2=\rho_i^2=\rho_i\sigma_i\lambda_i,
\label{squares}\\
e_i\sigma_j=\sigma_je_i, j\neq i,i-1, e_i\sigma_i=\sigma_ie_{i+1},
\sigma_ie_i=e_{i+1}\sigma_i, \label{com-e-sigma}\\
e_i\lambda_j=\lambda_je_i, j\neq i, i-1,
e_{i+1}\lambda_i=\lambda_i,
e_i\lambda_i=\lambda_ie_{i+1}=\lambda_ie_i=e_ie_{i+1},\label{two-zeros-lambda}\\
e_i\rho_j=\rho_je_i, j\neq i, i-1, \rho_ie_{i+1}=\rho_i,
\rho_ie_i=e_{i+1}\rho_i=e_i\rho_i=e_ie_{i+1}.
\label{two-zeros-rho}
\end{gather}

{\bf Remark 1.} It follows from \eqref{df-e} that $S$ is generated
by $\lambda_i$-s, $\rho_i$-s and $\sigma_i$-s only since the
defining relations can be readily rewritten without $e_i$-s.
However, it is convenient for us to include $e_i$-s to the
generating set and to the relations, because the products of
$e_i$-s will appear afterwards in the canonical words we will
introduce.

{\bf Remark 2.}  We would like to emphasize that the proposed set
of relations does not pretend to being irreducible. For example,
the relations $e_i^2=e_i$ from~\eqref{squares} follow from the
relations~\eqref{inverse}. However, we keep these and some other
redundant relations with the purpose of making the subsequent text
more transparent and readable.

In view of the relations~\eqref{c1}, \eqref{c2}, \eqref{c3} the
submonoid of $S$, generated by all $\sigma_i$-s,  is isomorphic to
the full symmetric group $\S_n$. From now on, identify this
submonoid with $\S_n$.

Let $\pi\in \S_n$ and $\alpha\in S$. Set
$\alpha^{\pi}=\pi^{-1}\alpha\pi$. Obviously, the map
$\varphi_{\pi}:\alpha\mapsto\alpha^{\pi}$ is an automorphism of
$S$, and $\pi\mapsto \varphi^{\pi}$ is an action of $\S_n$ on $S$.
Call this action {\em the action by inner automorphisms}.

Let $1\leq i<j\leq n$. Set
\begin{equation*}\sigma_{i,j}=\left\lbrace\begin{array}{l}\sigma_i, \text{ if
}
j=i+1;\\
\sigma_i\sigma_{i+1}\dots\sigma_{j-2}\sigma_{j-1}\sigma_{j-2}\dots
\sigma_i, \text{ if }j>i+1.\end{array}\right.
\end{equation*}

For $1\leq j<i\leq n$ we set $\sigma_{j,i}=\sigma_{i,j}$. Notice
that $\sigma_{i,j}^2=e$ for all acceptable $i,j$.

\begin{lemma}\label{lem:stabilizer}
Let $1\leq i\leq n-1$. Consider the action of $\S_n$ on $S$ by
inner automorphisms.
\begin{enumerate}
\item If $2\leq i\leq n-2$, then the elements $\sigma_1,\dots, \sigma_{i-2}, \sigma_{i+2},\dots,
\sigma_n$ and $\sigma_{i-1, i+2}$ stabilize both $\lambda_i$ and
$\rho_i$.
\item The elements $\sigma_3$, $\dots$, $\sigma_n$ stabilize both $\lambda_1$
and $\rho_1$, the elements $\sigma_1$, $\dots$, $\sigma_{n-2}$
stabilize both $\lambda_n$ and $\rho_n$
\item The elements $\sigma_1,\dots, \sigma_{i-2}, \sigma_{i+1},\dots,
\sigma_n$ and $\sigma_{i-1, i+1}$ stabilize $e_i$, $2\leq i\leq
n-1$.
\item The elements $\sigma_2$, $\dots$, $\sigma_{n-1}$ stabilize
$e_1$, the elements $\sigma_1$, $\dots$, $\sigma_{n-2}$ stabilize
$e_n$.
\end{enumerate}
\end{lemma}
\begin{proof}
To prove the first claim, in view of~\eqref{slr3}, we have only to
show that $\sigma_{i-1, i+2}$ stabilizes $\lambda_i$ and $\rho_i$,
$2\leq i\leq n-2$. Applying subsequently~\eqref{slr1}, \eqref{c1},
\eqref{slr3}, \eqref{c1}, \eqref{slr1}, \eqref{c1}, we obtain
\begin{multline*}
\sigma_{i-1,i+2}\lambda_i\sigma_{i-1,i+2}=
\sigma_{i-1}\sigma_i\sigma_{i+1}\sigma_i\sigma_{i-1}\lambda_i\sigma_{i-1}\sigma_i\sigma_{i+1}\sigma_i\sigma_{i-1}=\\
\sigma_{i-1}\sigma_i\sigma_{i+1}\lambda_{i-1}\sigma_{i+1}\sigma_i\sigma_{i-1}=
\sigma_{i-1}\sigma_i\lambda_{i-1}\sigma_i\sigma_{i-1}=\lambda_i,
\end{multline*}
as required. For $\rho_i$ arguments are similar.

To prove the third claim we let $2\leq i\leq n-1$ and show that
$\sigma_{i-1,i+1}$ stabilizes $e_i$. Indeed,
using~\eqref{com-e-sigma} we compute
\begin{multline*}
\sigma_{i-1}\sigma_i\sigma_{i-1}e_i\sigma_{i-1}\sigma_i\sigma_{i-1}=
\sigma_{i-1}\sigma_ie_{i-1}\sigma_{i}\sigma_{i-1}=\sigma_{i-1}e_{i-1}\sigma_{i}\sigma_{i}\sigma_{i-1}=\\
\sigma_{i-1}e_{i-1}\sigma_{i-1}=e_i,
\end{multline*}
as required.

The remaining two claims are proved similarly, and we leave the
details to the reader.
\end{proof}

To proceed, we need to introduce some more notation. Let $1\leq
p,q\leq n$ and $p\neq q$. For any $\pi\in \S_n$ such that
$\pi(1)=p$ and $\pi(2)=q$ set
\begin{equation}\label{def-lambda}
\lambda_{p,q}=\pi^{-1}\lambda_1\pi, \,\,
\rho_{p,q}=\pi^{-1}\rho_1\pi.
\end{equation}
In view of Lemma~\ref{lem:stabilizer} this definition is correct,
i.e. independent on the choice of $\pi\in \S_n$ such that
$\pi(1)=p$ and $\pi(2)=q$. Moreover, it can be easily verified
that
$$
\lambda_{i,i+1}=\lambda_i, \,\, \rho_{i,i+1}=\rho_i
$$
for all $1\leq i\leq n-1$. Indeed, for $i=1$ this is trivial. Let
$i\geq 2$. Then we apply~\eqref{slr1} and~\eqref{c1} $(i-1)$ times
and obtain
$$
(\sigma_{i-1}\sigma_i)\cdots (\sigma_2\sigma_3)
(\sigma_1\sigma_2)\lambda_1(\sigma_2\sigma_1)(\sigma_3\sigma_2)\cdots
(\sigma_i\sigma_{i-1})=\lambda_i.
$$
Besides, the element $(\sigma_2\sigma_1)(\sigma_3\sigma_2)\cdots
(\sigma_i\sigma_{i-1})$ maps $1$ to $i$ and $2$ to $i+1$
respectively.

\begin{lemma}\label{lem:lr}
Let $\pi\in\S_n$ be such that $\pi(p)=s$ and $\pi(q)=t$. Then
$\pi^{-1}\lambda_{p,q}\pi=\lambda_{s,t}$ and
$\pi^{-1}\rho_{p,q}\pi=\rho_{s,t}$.
\end{lemma}

\begin{proof}
We prove only the first equality, the second one being proved
similarly. Firstly, we note that every element $\alpha\in\S_n$
such that $\alpha(s)=s$ and $\alpha(t)=t$ stabilizes
$\lambda_{s,t}$. This follows from the definition of
$\lambda_{s,t}$ and Lemma~\ref{lem:stabilizer}. Further, consider
$\gamma$ and $\delta$ from $\S_n$ such that $\gamma(1)=s$,
$\gamma(2)=t$, $\delta(1)=p$, $\delta(2)=q$. Then
\begin{equation*}
\delta^{-1}\gamma\lambda_{p,q}\gamma^{-1}\delta=\delta^{-1}\lambda_1\delta=\lambda_{s,t},
\end{equation*}
which yields the required statement.
\end{proof}

\begin{lemma}\label{lem:e}
Let $\pi\in\S_n$ be such that $\pi(p)=s$. Then
$\pi^{-1}e_{p}\pi=e_{s}$.
\end{lemma}

\begin{proof}
Let $\alpha\in\S_n$ be such that $\alpha(1)=p$. Then
$\alpha^{-1}e_1\alpha=e_p$ by~\eqref{df-e},~\eqref{com-e-sigma}
and Lemma~\ref{lem:stabilizer}. Then we apply the arguments
similar to those from the proof of the previous lemma.
\end{proof}

In the following proposition we collect the relations satisfied by
the products of elements $\lambda_{p,q}$, $\rho_{p,q}$,
$\sigma_{p,q}$ by $e_i$.

\begin{proposition} \label{prop:rel_lamb} The following relations hold for all
admissible and pairwise distinct $p,q,k$:
\begin{gather}
e_p^2=e_p, e_p e_q=e_q e_p, \label{r1-00}\\
e_k\sigma_{p,q}=\sigma_{p,q}e_k,\,\,
e_p\sigma_{p,q}=\sigma_{p,q}e_q, \label{r1-0}\\
e_k\lambda_{p,q}=\lambda_{p,q}e_k,\,\,
e_q\lambda_{p,q}=\lambda_{p,q},\,\,
e_p\lambda_{p,q}=\lambda_{p,q}e_q=\lambda_{p,q}e_p=e_pe_q\label{r1-1}\\
e_k\rho_{p,q}=\rho_{p,q}e_k,\,\, \rho_{p,q}e_q=\rho_{p,q},\,\,
\rho_{p,q}e_p=e_q\rho_{p,q}=e_p\rho_{p,q}=e_pe_q \label{r1-2}
\end{gather}
\end{proposition}

\begin{proof} Consider the map $':S\to S$ defined by
$\sigma_i'=\sigma_i$, $\lambda_i'=\rho_i$, $\rho_i'=\lambda_i$,
$e_i'=e_i$. In view of the defining relations this map uniquely
extends to an involution on $S$ which we will also denote by $'$.

The relations~\eqref{r1-00} follow from~\eqref{squares} and
Lemma~\ref{lem:e}.

The relations~\eqref{r1-0} follow from \eqref{com-e-sigma} and
Lemma~\ref{lem:e}.

The relations \eqref{r1-1} follow from \eqref{two-zeros-lambda},
applying Lemma \ref{lem:e}.

Finally, \eqref{r1-2} follows from \eqref{r1-1}, using $'$.
\end{proof}

\begin{proposition}\label{pr:crucial}
For all pairwise distinct $p,q,k,l$,
\begin{gather}
\lambda_{p,q}^2=\rho_{p,q}^2=\lambda_{p,q}\lambda_{q,p}=\rho_{p,q}\rho_{q,p}=e_p e_q,\label{r1}\\
\lambda_{k,l}\lambda_{p,q}=\lambda_{p,q}\lambda_{k,l},\,\,
\rho_{k,l}\rho_{p,q}=\rho_{p,q}\rho_{k,l},\,\,
\lambda_{k,l}\rho_{p,q}=\rho_{p,q}\lambda_{k,l}, \label{r2}\\
\lambda_{k,q}\lambda_{k,l}=\lambda_{k,l}\lambda_{k,q}=\lambda_{k,q}\lambda_{q,l},
\rho_{k,l}\rho_{k,q}=\rho_{k,q}\rho_{k,l}=\rho_{l,q}\rho_{k,l}, \label{r3}\\
\lambda_{k,l}\lambda_{p,k}=e_l\lambda_{p,k},\,\,
\rho_{p,k}\rho_{k,l}=e_l\rho_{p,k},
\label{r4}\\
\lambda_{k,l}\lambda_{p,l}=e_k\lambda_{p,l}, \,\,
\rho_{p,l}\rho_{k,l}=e_k\rho_{p,l},\label{r5}\\
\lambda_{k,l}\rho_{l,k}=e_l\sigma_{k,l},\,\,
 \lambda_{k,l}\rho_{k,l}=e_l,\,\, \rho_{k,l}\lambda_{l,k}=e_ke_l, \label{r6}\\
\lambda_{k,l}\rho_{k,q}=\rho_{k,q}\lambda_{k,l}=\rho_{k,q}\lambda_{k,q}e_l\sigma_{q,l},\,\,
 \rho_{k,l}\lambda_{p,k}=\lambda_{p,k}e_l,\label{r7}\\
\lambda_{k,l}\rho_{p,k}=\rho_{p,k}\lambda_{p,k}e_l\sigma_{k,l},\,\,\lambda_{p,k}\rho_{k,l}=
\sigma_{k,l}e_l\rho_{p,k}\lambda_{p,k}\,\,
\rho_{p,k}\lambda_{k,l}=e_l\rho_{p,k},\label{r8}\\
\lambda_{k,l}\rho_{p,l}=\sigma_{p,l}e_p\rho_{k,l}\lambda_{k,l}\sigma_{p,l},\label{r10}\\
\rho_{k,l}\lambda_{k,l}=\rho_{l,k}\lambda_{l,k}\label{rr6}.
\end{gather}
\end{proposition}

\begin{proof}
To prove \eqref{r1}, in view of \eqref{squares} and applying
Lemmas \ref{lem:e}, \ref{lem:lr} and $'$, it is enough to check
that $\lambda_1\lambda_{2,1}=e_1e_2$. Indeed, using \eqref{slr2}
and \eqref{squares}, we have
\begin{equation*}
\lambda_1\lambda_{2,1}=\lambda_1\sigma_1\lambda_1\sigma_1=
\lambda_1^2=e_1e_2.
\end{equation*}

The relations \eqref{r2} follow from \eqref{lr1} and Lemma
\ref{lem:lr}.

To prove \eqref{r3}, in view of Lemma \ref{lem:lr} and applying
$'$, it is enough to prove that
$\lambda_1\lambda_{1,3}=\lambda_{1,3}\lambda_{1}=\lambda_1\lambda_2$.
Indeed, applying~\eqref{slr2},~\eqref{ghost} and
then~\eqref{ghost}, \eqref{slr2}, \eqref{ghost}, \eqref{slr2},
\eqref{slr1}, \eqref{slr1} we compute
\begin{equation*}
\lambda_1\lambda_{1,3}=\lambda_1\sigma_1\lambda_{2}\sigma_1=\lambda_1\lambda_{2}\sigma_1=
\lambda_1\lambda_{2};
\end{equation*}
\begin{multline*}
\lambda_1\lambda_{2}=\sigma_2\lambda_1\lambda_{2}=\sigma_2\lambda_1\lambda_{2}\sigma_2=
\sigma_2\lambda_1\lambda_{2}\sigma_1\sigma_2=\sigma_2\lambda_1\sigma_1\lambda_2\sigma_1\sigma_2=\\
\sigma_2\lambda_1\sigma_2\lambda_1\sigma_2\sigma_2=\sigma_1\lambda_2\sigma_1\lambda_1=
\lambda_{1,3}\lambda_1.
\end{multline*}

To prove~\eqref{r4} it is enough to show that
$\lambda_2\lambda_1=\lambda_1e_3$, which holds by~\eqref{ghost1}.

To prove~\eqref{r5} it is enough to check that
$\lambda_2\lambda_{1,3}=e_2\lambda_{1,3}$. Conjugating both sides
with $\sigma_2$ we obtain the equivalent equality
$\sigma_2\lambda_2\lambda_1=e_3\lambda_1$, which holds
by~\eqref{ghost1}.

The first equality of~\eqref{r6} follows from
$\lambda_2\rho_{3,2}=\lambda_2\rho_2\sigma_2=e_3\sigma_2$ and
Lemmas~\ref{lem:lr} and~\ref{lem:e}. The second and the third
equalities follow from the same lemmas in view of~\eqref{squares}.

In order to prove the first relation of \eqref{r7} it is again
enough to verify the relation
$\rho_{1,3}\lambda_1=\lambda_1\rho_{1,3}=\rho_{1,3}\lambda_{1,3}e_2\sigma_2$.
We compute
\begin{multline*}
\lambda_1\rho_{1,3}=\lambda_1\sigma_1\rho_2\sigma_1=\lambda_1\rho_2\sigma_1=
\sigma_2\rho_1\lambda_1e_3\sigma_1=\sigma_2\rho_1\lambda_1e_3=\sigma_2\rho_1\lambda_1
\sigma_2e_2\sigma_2=\\\sigma_2\rho_1\sigma_2\sigma_2\lambda_1\sigma_2e_2\sigma_2=
\sigma_1\rho_2\sigma_1\sigma_1\lambda_2\sigma_1e_2\sigma_2=\rho_{1,3}\lambda_{1,3}e_2\sigma_2
\end{multline*}
and
\begin{equation*}
\lambda_1\rho_{1,3}=\lambda_1\rho_2\sigma_1=\sigma_1\rho_2\sigma_1\lambda_1\sigma_1=
\rho_{1,3}\lambda_1\sigma_1=\rho_{1,3}\lambda_1.
\end{equation*}

To prove the second relation of \eqref{r7} it suffices  to
establish only that $\rho_2\lambda_1=\lambda_1e_3$, which is a
direct consequence of \eqref{neighb-rho-lambda}. Further, applying
$'$ to this relation we obtain the third relation of \eqref{r7}.

To prove the first relation of \eqref{r8} we verify only that
$\lambda_2\rho_1=\rho_1\lambda_1e_3\sigma_2$, which is a direct
consequence of \eqref{neighb-lambda-rho}. The second relation of
\eqref{r8} is a consequence of the first one using $'$.

 Further,
let us prove \eqref{r10}. It is enough to establish that
$\lambda_1\rho_{3,2}=\sigma_{2}e_3\rho_{1}\lambda_1\sigma_2$. This
equality is equivalent to
$\lambda_1\rho_2=\sigma_2e_3\rho_1\lambda_1$, which, in turn,
follows from \eqref{neighb-lambda-rho}.

 Finally, the
relation~\eqref{rr6} follows from
$\sigma_1\rho_1\sigma_1\sigma_1\lambda_1\sigma_1=\rho_1\lambda_1$
and Lemma~\ref{lem:lr}. The proof is complete.
\end{proof}

\section{Rewriting technique and canonical words in the monoid
$S$}\label{s4}

In this section we are going to develop some rewriting technique
in order to show that any element of the monoid $S$ can be
represented by some "reduced word". This will imply that the
cardinality of $S$ is not bigger then the cardinality of the set
of all reduced words.

We start from the following observation.

\begin{lemma}\label{lem:rewr1}
Every element of $S$ can be written as a product of the form
$\alpha_1\dots \alpha_k\beta$, where $k\geq 0$, each $\alpha_i$ is
equal to some $\lambda_{p,q}$ or $\rho_{p,q}$ and $\beta\in\S_n$.
\end{lemma}
\begin{proof}
Let $\alpha\in S$. Since every $\lambda_i$ and $\rho_i$ belongs to
some orbit of $\lambda_1$ and $\rho_1$ respectively with respect
to the action of $\S_n$ by inner automorphisms,  $S$ can be
generated by $\S_n$, $\lambda_1$ and $\rho_1$. Therefore, $\alpha$
can be written in the form
$$
\alpha=\pi_1\gamma_1\pi_2\gamma_2\dots\pi_k\gamma_k\pi_{k+1},
$$
where $k\geq 0$, $\pi_i\in\S_n$ for all $i$ and each $\gamma_i$
equals either $\lambda_1$ or $\rho_1$. In view
of~\eqref{def-lambda} we can rewrite the expression for $\alpha$
as follows:
\begin{multline*}
\alpha=\pi_1\gamma_1\pi_1^{-1}(\pi_1\pi_2)\gamma_2(\pi_1\pi_2)^{-1}\dots
(\pi_1\dots\pi_k)\gamma_k(\pi_1\dots\pi_k)^{-1}\cdot\\\cdot(\pi_1\dots\pi_k\pi_{k+1})=
\gamma^1_{p_1,q_1}\dots\gamma^k_{p_k,q_k}\beta,
\end{multline*}
where $p_i=(\pi_1...\pi_{i})^{-1}(1)$,
$q_i=(\pi_1...\pi_{i})^{-1}(2)$ and
\begin{equation*}
\gamma^i_{p_i,q_i}=\left\lbrace\begin{array}{l}\lambda_{p_i,q_i},
\text{ if } \gamma_i=\lambda_1,\\
\rho_{p_i,q_i},\text{ if }\gamma_i=\rho_1,\end{array}\right.
\end{equation*}

$1\leq i\leq k$, and $\beta=\pi_1\dots\pi_{k+1}$.
\end{proof}

\begin{lemma}\label{lem:rewr2} Every element of $S$ can be written as a product of the
form
\begin{equation}\label{eq:product}
\alpha_{1}\dots\alpha_k\beta_1\dots\beta_l\pi E,
\end{equation}
where $k,l\geq 0$, each $\alpha_i$ equals some $\rho_{p,q}$, each
$\beta_i$ equals some $\lambda_{p,q}$, $\pi\in\S_n$ and $E=e$ or
$E$ is a product of several $e_p$-s.
\end{lemma}
\begin{proof}
Let $\alpha\in S$. It follows from Lemma~\ref{lem:rewr1} that we
can express $\alpha$ as a product $\alpha_1\dots \alpha_k\pi$,
where $k\geq 0$, each $\alpha_i$ is equal to some $\lambda_{p,q}$
or $\rho_{p,q}$ and $\pi\in\S_n$. Suppose that
$\alpha_i=\lambda_{p,q}$ and $\alpha_{i+1}=\rho_{k,l}$ for some
$i$. If the sets $\{p,q\}$ and $\{k,l\}$ are disjoint we have
$\alpha_i\alpha_{i+1}=\alpha_{i+1}\alpha_i$ by~\eqref{r2}. If the
sets $\{p,q\}$ and $\{k,l\}$ are not disjoint we apply the
appropriate relation of~\eqref{r6}-\eqref{r10}. As a result we
obtain an expression for $\alpha$ containing less subwords of the
form $\lambda_{i,j}\rho_{s,t}$.

However, after such a rewriting some $e_i$-s and $\sigma_{s,t}$-s
might appear in the expression for $\alpha$. If some $e_i$-s
appear, using~\eqref{r1-0}, \eqref{r1-1} and~\eqref{r1-2} we
rewrite our expression such that it has the occurrence of $e_j$ at
the rightmost position, while the number of subwords of the form
$\lambda_{i,j}\rho_{s,t}$ remains the same. If some
$\sigma_{s,t}$-s appear, using the action of $\S_n$ on $S$ by
inner automorphisms, we can, similarly to as this is done in the
proof of Lemma~\ref{lem:rewr1}, rewrite it such that the group
element occurs to the right to all occurrences of
$\lambda_{i,j}$-s and $\rho_{s,t}$-s. As the mentioned rewriting
does not affect the number of the subwords of the form
$\lambda_{i,j}\rho_{s,t}$, the statement of the lemma follows by
induction on the number of subwords of the form
$\lambda_{i,j}\rho_{s,t}$ in the initial expression for $\alpha$.
\end{proof}

We can even strengthen the previous statement.

\begin{lemma}\label{lem:rewr3} Every element of $S$ can be written
as a product of the form~\eqref{eq:product} such that the
conditions of Lemma~\ref{lem:rewr2} are satisfied and, moreover,
the following conditions are also satisfied:
\begin{enumerate}
\item If $\alpha_i=\rho_{p,q}$ and $\alpha_{j}=\rho_{k,l}$ then
either $\{p,q\}\cap\{k,l\}=\varnothing$ or
$\{p,q\}\cap\{k,l\}=\{p\}=\{k\}$, so that $\alpha_i=\rho_{p,q}$
and $\alpha_{j}=\rho_{p,l}$. In particular, $\alpha_i$ and
$\alpha_{j}$ commute.
\item If $\beta_i=\lambda_{p,q}$ and $\beta_{j}=\lambda_{k,l}$ then
either $\{p,q\}\cap\{k,l\}=\varnothing$ or
$\{p,q\}\cap\{k,l\}=\{p\}=\{k\}$, so that $\beta_i=\lambda_{p,q}$
and $\beta_{j}=\lambda_{p,l}$. In particular, $\beta_i$ and
$\beta_{j}$ commute.
\end{enumerate}
\end{lemma}
\begin{proof}
We will prove the statement on $\alpha_i$-s only, the second
statement being proved analogously.

Notice that it is enough to prove the statement for the case
$j=i+1$. Indeed, if $\alpha_i$ commutes with $\alpha_{i+1}$ for
every $i$ then we can rearrange the factors of $\alpha$ such that
$\alpha_j$ follows $\alpha_i$ for any $i,j$.

Apply induction on the number of factors of the form $\rho_{p,q}$
in the expression~\eqref{eq:product}. If this number is zero or
one, the statement is obvious. Suppose that $\alpha_1=\rho_{p,q}$
and $\alpha_{2}=\rho_{k,l}$, where $p,q,l$ are pairwise distinct,
and that $\{p,q\}\cap\{k,l\}\neq\varnothing$. Consider six
possible cases:

\begin{enumerate}[a)]
\item \label{a1} if $\alpha_i=\rho_{p,q}, \alpha_{i+1}=\rho_{p,q}$,
we apply~\eqref{r1},
\item  if \label{a2} $\alpha_i=\rho_{p,q}, \alpha_{i+1}=\rho_{q,p}$, we apply~\eqref{r1},
\item if \label{a3} $\alpha_i=\rho_{p,q}, \alpha_{i+1}=\rho_{p,l}$,
then
$\alpha_i\alpha_{i+1}=\rho_{p,q}\rho_{p,l}=\rho_{p,l}\rho_{p,q}$
by~\eqref{r3},
\item if \label{a4} $\alpha_i=\rho_{p,q}, \alpha_{i+1}=\rho_{l,p}$, then  the product $\alpha_i\alpha_{i+1}$ equals the
product from \ref{a3}) by~\eqref{r3},
\item if \label{a5} $\alpha_i=\rho_{p,q}, \alpha_{i+1}=\rho_{l,q}$,
we apply~\eqref{r5},
\item if \label{a6} $\alpha_i=\rho_{p,q}, \alpha_{i+1}=\rho_{q,l}$,
we apply~\eqref{r4}.
\end{enumerate}

In the cases~\ref{a1}), \ref{a2}), \ref{a5}), \ref{a6}) we obtain
an expression for the initial element containing less entries of
factors of the form $\rho_{p,q}$ and apply the inductive
hypothesis. In the cases~\ref{a3}) and \ref{a4}) we have that
$\alpha_i\alpha_{i+1}=\rho_{p,q}\rho_{p,k}$ for some pairwise
distinct $p,q,k$.

We proceed by considering the product $\alpha_2\alpha_3$ and so
on. Finally we either reach the last factor $\alpha_l$ with the
first claim satisfied for every possible $\alpha_i$ and
$\alpha_{i+1}$, or reduce the number of factors. In the latest
case we apply the induction.
\end{proof}

Let $p\in\{1,\dots, n\},$ $A\subseteq \{1,\dots, n\}$ and
$p\not\in A$. Set $R_{p, A}=\rho_{p,a_1}\cdots \rho_{p, a_s}$ and
$L_{p,A}=\lambda_{p,a_1}\cdots\lambda_{p,a_s}$, where
$A=\{a_1,\dots, a_s\}$. In view of~\eqref{r3} $R_{p,A}$ and
$L_{p,A}$ are well-defined.

\begin{corollary}\label{cor:commutation} Suppose $A_1\cap
A_2=\varnothing$, $p_1\neq p_2$, $p_1\not\in A_1$, $p_2\not\in
A_2$. Then $R_{p_1,A_1}R_{p_2,A_2}=R_{p_2,A_2}R_{p_1,A_1}$ and
$L_{p_1,A_1}L_{p_2,A_2}=L_{p_2,A_2}L_{p_1,A_1}$.
\end{corollary}

For a subset $M=\{m_1,\dots, m_s\}\subset \{1,\dots, n\}$ set
$E_M=e_{m_1}\cdots e_{m_s}$, which is well-defined in view
of~\eqref{r1-00}.

\begin{proposition}\label{prop:rewr4} Every element of $S$ can be written
as a product of the form
\begin{equation}\label{eq:expression}
R_{p_1, A_1}\dots R_{p_k, A_k} L_{q_1, B_1}\dots L_{q_l, B_l} E_M
\sigma,
\end{equation}
where $k,l\geq 0$, $p_1,\dots, p_k, q_1,\dots, q_l$ are pairwise
distinct and $A_1,\dots, A_k$, $B_1, \dots, B_l$ are pairwise
disjoint, $E=E_M$, $M\subseteq \{1,\dots, n\}$, $\sigma\in\S_n$.
Moreover, the following conditions are satisfied:
\begin{enumerate}[(i)]
\item \label{i}$p_i\not\in B_1\cup\dots \cup B_l$, $1\leq i\leq k$,

\item \label{ii}$q_i\not\in A_1\cup\dots\cup A_k$, $1\leq i\leq l$,

\item \label{iii} $M$ is disjoint with $\bigl(\cup_{i=1}^{k}(\{p_i\}\cup
A_i)\bigr)\cup\bigl(\cup_{i=1}^l (\{q_i\}\cup B_i)\bigr)$.
\end{enumerate}
\end{proposition}
\begin{proof}
Let $\alpha\in S$ be presented in the form \eqref{eq:product},
such that the conditions of Lemma~\ref{lem:rewr3} are satisfied.
The relations~\eqref{r1-0} imply that we can move from the
expression \eqref{eq:product} to
\begin{equation}\label{eq:product_new}
\alpha=\alpha_1\cdots\alpha_k\beta_1\cdots \beta_l E'\pi',
\end{equation}
where $E'$ is the product of some $e_i$-s and $\pi'\in\S_n$. It
follows from Lemma~\ref{lem:rewr3} that we can rearrange the
factors in the product $\alpha_1\cdots \alpha_k$ and obtain the
expression
$$
\alpha_1\cdots \alpha_k=R_{p_1, A_1}\cdots R_{p_s, A_s}
$$
for certain pairwise distinct $p_1,\dots, p_s$ and pairwise
disjoint $A_1,\dots, A_s$. Similarly, we can rearrange the factors
in the product $\beta_1\cdots \beta_l$ such that it equals
$L_{q_1, B_1}\cdots L_{q_t, B_t}$ for certain pairwise distinct
$q_1,\dots, q_t$ and pairwise disjoint $B_1,\dots, B_t$. It
follows that we can achieve the expression of the form
~\eqref{eq:expression} for $\alpha$.

Suppose $p_i\in B_j$, for some $1\leq i\leq k$ and $1\leq j\leq
l$. Rearranging, if necessary, the factors in $\alpha$, we obtain
an expression for $\alpha$ containing the factor $R_{p_i,
A_i}L_{q_j,B_j}$. Then rearrange the factors constituting
$L_{q_j,B_j}$ in such a way that the obtained factorization of
$\alpha$ contains the factor $R_{p_i, A_i}\lambda_{q_j,p_i}$. In
view of~\eqref{r7}
$\rho_{p_i,a}\lambda_{q_j,p_i}=\lambda_{q_j,p_i}e_a$. Applying
this equality several times we obtain
\begin{equation*}
R_{p_i, A_i}\lambda_{q_j,p_i}=\lambda_{q_j,p_i}\prod_{a\in
A_i}e_a.
\end{equation*}

Applying~\eqref{r1-1} to the obtained expression for $\alpha$ we
move all $e_a$-s to the right of all $L_{q_j,B_j}$-s. The
resulting expression for $\alpha$ will be of the
form~\eqref{eq:expression},  but without the factor $R_{p_i,A_i}$,
moreover possibly without some $L_{q_j,B_j}$-s and with a new $E$
(containing more $e_i$-s). What we have reached is that the number
of $p_i$-s contained in some $B_j$-s in the renewed expression for
$\alpha$ is decreased by one. Applying the described rewriting
several times we obtain an expression for $\alpha$ such that the
condition~{\em(\ref{i})} is satisfied.

Applying analogous manipulations the expression for $\alpha$ can
be rewritten such that the condition~{\em(\ref{ii})} is also
satisfied.

Now we can assume that $\alpha$ is written in the
form~\eqref{eq:expression} and the conditions~{\em(\ref{i})}
and~{\em(\ref{ii})} are satisfied. Suppose there is $a\in M$ such
that $a=q_i$ or $a\in B_i$ for certain $i$. Rewrite the expression
for $\alpha$ such that it contains the factor $L_{q_i,B_i}e_{a}$
and apply~\eqref{r1-1} several times. We will obtain
$$L_{q_i,B_i}e_{a}=e_{q_i}\prod_{b\in B_i}e_b.$$
This and inductive arguments on the number of factors of the form
$L_{q_i,B_i}$ in~\eqref{eq:expression} show that $\alpha$ can be
rewritten such that in the given expression for $\alpha$ the set
$M$ is disjoint with $\cup_{i=1}^l (\{q_i\}\cup B_i)$.

Let us continue the rewriting of the expression for $\alpha$.
Suppose there is $a\in M$ such that $a=p_i$ or $a\in A_i$ for
certain $i$. We can assume that $e_a$ is the first factor of $E$.
In view of the first relation of~\eqref{r1-1} $e_a$ commutes with
every $L_{q_i, B_i}$. Hence we rewrite the expression such that
$e_a$ is located between the group of the factors $R_{p_i, A_i}$-s
and the group of the factors $L_{q_i, B_i}$-s of our
expression~\eqref{eq:expression}.  Moreover, we can assume that
this expression contains the factor $R_{p_i,A_i}e_a$. Similarly to
as it was done previously we rearrange this factor and obtain
$$e_{a}R_{p_i,A_i}=e_{p_i}\prod_{x\in A_i}e_x.$$

Thus the number of such $a\in M$ that $a=p_i$ or $a\in A_i$ for
certain $i$ has been decreased by one. The difficulty here is that
the current expression for $\alpha$ may be not of the
form~\eqref{eq:expression}. To reach the expression of the
required form we have to move the product $e_{p_i}\prod_{x\in
A_i}e_x$ to the position to the right of  all the $L_{q_i,
B_i}$-s. It is enough to show that such a movement is possible for
every $e_x$ with $x\in \{p_i\}\cup A_i$ and apply induction. If
$x\not\in \cup_{i=1}^l (\{q_i\}\cup B_i)$ then $e_x$ commutes with
each $L_{q_i, B_i}$, and the required movement is performed.
Otherwise, applying (possibly several times)~\eqref{r1-1} we
obtain $$e_xL_{q_i, B_i}=e_{q_i}\prod_{b\in B_i}e_b.$$

Since the sets $\{q_j\}\cup B_j$, $1\leq j\leq l$, are pairwise
disjoint it follows that every $y\in \{q_i\}\cup B_i$ does not
belong to any of the sets $\{q_j\}\cup B_j$, $1\leq j\leq l$,
$j\neq i$. Therefore, $e_y$ commutes with all $L_{q_j, B_j}$-s,
$j\neq i$, by~\eqref{r1-1}. This completes the proof.
\end{proof}

Denote
$T=\bigl(\cup_{i=1}^k\{p_i\}\bigr)\cup\bigl(\cup_{i=1}^l\{q_i\}\bigr)$,
$s=|T|$. Enumerate the elements of $T$ in some way, suppose
$T=\{t_1,\dots, t_s\}$.

 Define the sets $C_i$ and $R^{t_i}_{C_i}$, $1\leq i\leq s$, in the
following way. Let $1\leq i\leq s$.
\begin{itemize}
\item If $t_i=p_j$ for some $j$ we set
$C_i=A_j\cup\{p_j\}$ and $R^{t_i}_{C_i}=R_{p_j, A_j}$.
\item If $t_i\not\in\cup_{j=1}^k\{p_j\}$ we set
$C_i=\{t_i\}$ and $R^{t_i}_{C_i}=e$.
\end{itemize}

The above defined sets $C_1$, $\dots$, $C_s$ are pairwise
disjoint, their union coincides with $T\cup (\cup_{i=1}^k A_i)$.
Moreover, $t_i\in C_i$ for every possible $i$.

Similarly, define the sets $D_i$,  $L^{t_i}_{D_i}$, $1\leq i\leq
s$:
\begin{itemize}
\item If $t_i=q_j$ for some $j$ we set
$D_i=B_j\cup\{q_j\}$ and $L^{t_i}_{D_i}=L_{q_j, B_j}$.
\item If $t_i\not\in\cup_{j=1}^k\{q_j\}$ we set
$D_i=\{t_i\}$ and $L^{t_i}_{D_i}=e$.
\end{itemize}

The above defined sets $D_1$, $\dots$, $D_s$ are pairwise
disjoint, their union coincides with $T\cup \bigl(\cup_{i=1}^k
B_i\bigr)$. Moreover, $t_i\in D_i$ for every possible $i$.

\begin{corollary}\label{cor:can_form} Every element of $S$ can be presented in the
form
\begin{equation}\label{eq:canonical}
R^{t_1}_{C_1}\cdots R^{t_s}_{C_s}L^{t_1}_{D_1}\cdots
L^{t_s}_{D_s}E_M\sigma,
\end{equation}
where $C_1, \dots C_s$ are pairwise disjoint, $D_1, \dots D_s$ are
pairwise disjoint, $t_i\in C_i\cap D_i$, $1\leq i\leq s$,  $M$ is
disjoint with $\bigl(\cup_{i=1}^s (C_i\cup D_i)\bigr)$ and
$\sigma\in\S_n$.
\end{corollary}

\begin{proof} Follows from Proposition~\ref{prop:rewr4}.
\end{proof}

Let $B$ be a subset of $\{1,\dots, n\}$. Denote by $\S_B$ the
subgroup of $\S_n$ generated by all $\sigma_{i,j}$ with $i,j\in
B$. Obviously, $\S_B$ is isomorphic to $\S_{|B|}$.

Call an expression of the form~\eqref{eq:canonical} such that the
conditions of Corollary~\ref{cor:can_form} are satisfied a {\em
canonical word}.

Let $F=M\cup (\cup_{i=1}^s(C_i\setminus D_i))$ and
$G=\S_{D_1}\oplus\cdots\oplus \S_{D_s}\oplus \S_F$. The group $G$
depends on $\{D_1,\dots, D_s\}, F$, but we do not indicate this
into the notation just not to overload it.

Call two canonical words
\begin{equation*}
R^{t_1}_{C_1}\cdots R^{t_s}_{C_s}L^{t_1}_{D_1}\cdots
L^{t_s}_{D_s}E_{M_1}\sigma_1 {\text{ and }}
\end{equation*}
\begin{equation*}
R^{t'_1}_{C'_1}\cdots R^{t'_s}_{C'_s}L^{t'_1}_{D'_1}\cdots
L^{t'_s}_{D'_s}E_{M_2}\sigma_2
\end{equation*}
{\em equivalent} provided that there is a permutation
$\tau\in\S_s$ such that $C_i=C'_{\tau(i)}$, $D_i=D'_{\tau(i)}$,
$1\leq i\leq s$, $M_1=M_2$ and $\sigma_1\sigma_2^{-1}\in G$.

\begin{proposition}\label{prop:canonical} If two canonical words are equivalent then
their values in $S$ are equal.
\end{proposition}

The proof will be derived from a series of the following lemmas.

\begin{lemma}\label{lem:aux}
For pairwise distinct $i,j,q$
\begin{equation}\label{eq:eat_lambda}
\lambda_{q,i}\lambda_{q,j}\sigma_{i,j}=\lambda_{q,i}\lambda_{q,j}
\end{equation}
Furthermore,
\begin{equation}\label{eq:one_more}
\lambda_{q,i}\lambda_{q,j}\sigma_{q,j}=\lambda_{q,i}\lambda_{q,j}.
\end{equation}
\end{lemma}

\begin{proof}
 Applying~\eqref{r3} and~\eqref{slr2} we compute
\begin{equation*}
\lambda_{q,i}\lambda_{q,j}\sigma_{i,j}=\lambda_{q,i}\lambda_{i,j}\sigma_{i,j}=
\lambda_{q,i}\pi^{-1}\lambda_1\sigma_1\pi=
\lambda_{q,i}\pi^{-1}\lambda_1\pi=
\lambda_{q,i}\lambda_{i,j}=\lambda_{q,i}\lambda_{q,j},
\end{equation*}
where $\pi$ is an element of $\S_n$ such that $\pi(1)=i$ and
$\pi(2)=j$. The relation~\eqref{eq:one_more} follows by the same
argument.
\end{proof}

\begin{lemma} Let $i\neq j$. Then
\begin{equation}\label{eq:esigma}
e_ie_j\sigma_{i,j}=e_ie_j.
\end{equation}
\end{lemma}
\begin{proof}
Applying consequently~\eqref{r6},
$\rho_{i,j}\sigma_{i,j}=\rho_{j,i}$ (which holds by
Lemma~\ref{lem:lr} and~\eqref{slr2}),~\eqref{r6} and~\eqref{r1-1}
we compute
\begin{equation*}
e_ie_j\sigma_{i,j}=e_i\lambda_{i,j}\rho_{i,j}\sigma_{i,j}=\lambda_{j,i}\rho_{j,i}\lambda_{i,j}\rho_{j,i}=\lambda_{j,i}e_ie_j\rho_{j,i}=e_ie_j.
\end{equation*}
\end{proof}

\begin{lemma}\label{lem:stabil_S} Let $\alpha=R^{t_1}_{C_1}\cdots R^{t_s}_{C_s}L^{t_1}_{D_1}\cdots
L^{t_s}_{D_s}E_{M}.$  Then $\alpha$ is stabilized by $G$ from the
right, that is $\alpha\sigma=\alpha$ for every $\sigma\in G$.
\end{lemma}
\begin{proof}

 Suppose first that $\alpha\in \S_{D_r}$, $1\leq r\leq s$. It is
 enough to consider only the case when $\alpha=\sigma_{i,j}$,
 $i,j\in D_r$. Since $E_M$ commutes with $L^{t_r}_{D_r}$ we are
 only to establish that $L^{t_r}_{D_r}\sigma_{i,j}=L^{t_r}_{D_r}$.
 But this equality follows from~\eqref{eq:eat_lambda} and~\eqref{eq:one_more}.

 Suppose now that $\alpha\in \S_{F}$. As in the previous
 paragraph, we  consider only the case when $\alpha=\sigma_{i,j}$,
 $i,j\in F$. To prove the statement it is enough to show that
 \begin{equation}\label{eq:important}
 \alpha\sigma_{i,j}=\alpha e_ie_j\sigma_{i,j}
\end{equation}
and apply~\eqref{eq:esigma}.

If $i\in M$ then $\alpha e_i=\alpha$ by the definition of $E_M$
and~\eqref{squares}.

Suppose $i\in C_r\setminus D_r$ for some $r$, $1\leq r\leq s$, and
show that $\alpha e_i=\alpha$. Firstly, $E_M e_i=e_iE_M$
by~\eqref{squares}. Further, $L_{D_j}^{t_j}e_i=e_iL_{D_j}^{t_j}$,
$1\leq j\leq s$, by the first relation of~\eqref{r1-1}. Finally,
$R_{C_r}^{t_r}e_i=R_{C_r}^{t_r}$ by the second relation
of~\eqref{r1-2}. This completes the proof.
\end{proof}

\begin{proof}[Proof of Proposition~\ref{prop:canonical}.]
Suppose  $R^{t_1}_{C_1}\cdots R^{t_s}_{C_s}L^{t_1}_{D_1}\cdots
L^{t_s}_{D_s}E_{M}$ is a canonical word and $x_i\in C_i\cap D_i$,
$1\leq i\leq s$. In view of Lemma~\ref{lem:stabilizer} and
Corollary~\ref{cor:commutation} it is enough to show that if we
replace $t_i$-s by $x_i$-s, the obtained canonical word has the
same value in $S$. Fix some index $i$. We can assume that the
initial canonical word has the factor $\rho_{t_i,
x_i}\lambda_{t_i,x_i}$. In view of~\eqref{rr6} this factor equals
$\rho_{x_i, t_i}\lambda_{x_i,t_i}$. It follows that
$R^{t_i}_{C_i}L^{t_i}_{D_i}=R^{x_i}_{C_i}L^{x_i}_{D_i}$, which
completes the proof.
\end{proof}

\section{Canonical form for the elements of $\PIP_n$}\label{s5}

We start this section from introducing the notation for certain
elements of the monoid $\PIP_n$. For distinct $x$ and $y$ of $X$
we set
\begin{equation*}
s_{x,y}=\bigl\{\{x,y'\},\{x',y\},\{t,t'\}_{t\in
X\setminus\{x,y\}}\bigr\},
\end{equation*}
\begin{equation*}
r_{x,y}=\bigl\{\{x,y,x'\},\{y'\},\{t,t'\}_{t\in
X\setminus\{x,y\}}\bigr\},
\end{equation*}
\begin{equation*}
l_{x,y}=\bigl\{\{x,x',y'\},\{y\},\{t,t'\}_{t\in
X\setminus\{x,y\}}\bigr\}~\mbox{and}
\end{equation*}
\begin{equation*}
\varepsilon_{x}=\bigl\{\{x\},\{x'\},\{t,t'\}_{t\in
X\setminus\{x\}}\bigr\}.
\end{equation*}

Furthermore, we set $s_i=s_{i,i+1}$, $r_i=r_{i,i+1}$ and
$l_i=l_{i,i+1}$ for $1\leq i\leq n-1$. The elements $s_1, \dots,
s_{n-1}$ generate the group of units of $\PIP_n$ which is
isomorphic to the symmetric group $\S_n$ and will be identified
with it.

We will use the following statement.

\begin{proposition}[\cite{KMal}]\label{pr:gen-pip}
Let $n\geq 3$. Then $\PIP_n$ is generated by $\Sym_n$, $r_1$ and
$l_1$.
\end{proposition}

\begin{proposition}\label{prop:rel_pip_n}
The map from $S$ to $\PIP_n$, sending $\sigma_i$ to $s_i$,
$\lambda_i$ to $l_i$ and $\rho_i$ to $r_i$, $1\leq i\leq n-1$,
extends to an epimorphism $\varphi:S\to \PIP_n$.
\end{proposition}

\begin{proof} Firstly we make sure that
$\varepsilon_i=l_{i-1}r_{i-1}$, $2\leq i\leq n$, and
$\varepsilon_1=s_1\varepsilon_2 s_1$. Then we verify that for the
elements $s_i$, $l_i$, $r_i$ and $\varepsilon_i$ all the relations
corresponding to the relations~\eqref{c1}-\eqref{two-zeros-rho}
hold. This and Proposition~\ref{pr:gen-pip} imply the needed
statement.
\end{proof}

\begin{figure}
\special{em:linewidth 0.4pt} \unitlength 0.80mm
\linethickness{1pt}
\begin{picture}(150.00,35.00)

\put(34.00,00.00){\makebox(0,0)[cc]{$\bullet$}}
\put(34.00,10.00){\makebox(0,0)[cc]{$\bullet$}}
\put(34.00,20.00){\makebox(0,0)[cc]{$\bullet$}}
\put(44.00,00.00){\makebox(0,0)[cc]{$\bullet$}}
\put(44.00,10.00){\makebox(0,0)[cc]{$\bullet$}}
\put(44.00,20.00){\makebox(0,0)[cc]{$\bullet$}}

\drawline(30.00,-04.00)(30.00,24.00)
\drawline(30.00,24.00)(48.00,24.00)
\drawline(48.00,24.00)(48.00,-04.00)
\drawline(48.00,-04.00)(30.00,-04.00)

\drawline(34.30,00.30)(44.30,00.30)

\drawline(32.30,08.30)(32.30,12.30)
\drawline(32.30,12.30)(36.30,12.30)
\drawline(36.30,12.30)(36.30,08.30)
\drawline(36.30,08.30)(32.30,08.30)

\drawline(32.30,18.30)(32.30,22.30)
\drawline(32.30,22.30)(46.30,22.30)
\drawline(46.30,22.30)(46.30,08.30)
\drawline(46.30,08.30)(42.80,08.30)
\drawline(42.30,08.30)(32.30,18.30)


\put(57.00,00.00){\makebox(0,0)[cc]{$\bullet$}}
\put(57.00,10.00){\makebox(0,0)[cc]{$\bullet$}}
\put(57.00,20.00){\makebox(0,0)[cc]{$\bullet$}}
\put(67.00,00.00){\makebox(0,0)[cc]{$\bullet$}}
\put(67.00,10.00){\makebox(0,0)[cc]{$\bullet$}}
\put(67.00,20.00){\makebox(0,0)[cc]{$\bullet$}}

\drawline(53.00,-04.00)(53.00,24.00)
\drawline(53.00,24.00)(71.00,24.00)
\drawline(71.00,24.00)(71.00,-04.00)
\drawline(71.00,-04.00)(53.00,-04.00)

\drawline(55.30,-02.30)(55.30,02.30)
\drawline(55.30,02.30)(65.30,12.30)
\drawline(65.30,12.30)(69.30,12.30)
\drawline(69.30,12.30)(69.30,-02.30)
\drawline(69.30,-02.30)(55.30,-02.30)

\drawline(57.30,20.30)(67.30,20.30)

\drawline(55.30,08.30)(55.30,12.30)
\drawline(55.30,12.30)(59.30,12.30)
\drawline(59.30,12.30)(59.30,08.30)
\drawline(59.30,08.30)(55.30,08.30)


\put(85.00,00.00){\makebox(0,0)[cc]{$\bullet$}}
\put(85.00,10.00){\makebox(0,0)[cc]{$\bullet$}}
\put(85.00,20.00){\makebox(0,0)[cc]{$\bullet$}}
\put(95.00,00.00){\makebox(0,0)[cc]{$\bullet$}}
\put(95.00,10.00){\makebox(0,0)[cc]{$\bullet$}}
\put(95.00,20.00){\makebox(0,0)[cc]{$\bullet$}}

\drawline(81.00,-04.00)(81.00,24.00)
\drawline(81.00,24.00)(99.00,24.00)
\drawline(99.00,24.00)(99.00,-04.00)
\drawline(99.00,-04.00)(81.00,-04.00)

\drawline(85.30,00.30)(95.30,00.30)
\drawline(85.30,10.30)(95.30,10.30)

\drawline(83.30,18.30)(83.30,22.30)
\drawline(83.30,22.30)(87.30,22.30)
\drawline(87.30,22.30)(87.30,18.30)
\drawline(87.30,18.30)(83.80,18.30)

\drawline(93.30,18.30)(93.30,22.30)
\drawline(93.30,22.30)(97.30,22.30)
\drawline(97.30,22.30)(97.30,18.30)
\drawline(97.30,18.30)(93.80,18.30)


\put(108.00,00.00){\makebox(0,0)[cc]{$\bullet$}}
\put(108.00,10.00){\makebox(0,0)[cc]{$\bullet$}}
\put(108.00,20.00){\makebox(0,0)[cc]{$\bullet$}}
\put(118.00,00.00){\makebox(0,0)[cc]{$\bullet$}}
\put(118.00,10.00){\makebox(0,0)[cc]{$\bullet$}}
\put(118.00,20.00){\makebox(0,0)[cc]{$\bullet$}}

\drawline(104.00,-04.00)(104.00,24.00)
\drawline(104.00,24.00)(122.00,24.00)
\drawline(122.00,24.00)(122.00,-04.00)
\drawline(122.00,-04.00)(104.00,-04.00)

\drawline(106.30,-02.30)(106.30,02.30)
\drawline(106.30,02.30)(116.30,12.30)
\drawline(116.30,12.30)(120.30,12.30)
\drawline(120.30,12.30)(120.30,-02.30)
\drawline(120.30,-02.30)(106.30,-02.30)

\drawline(108.30,20.30)(118.30,20.30)

\drawline(106.30,08.30)(106.30,12.30)
\drawline(106.30,12.30)(110.30,12.30)
\drawline(110.30,12.30)(110.30,08.30)
\drawline(110.30,08.30)(106.30,08.30)


\drawline(44.30,00.30)(47.30,00.30)
\drawline(49.30,00.30)(52.30,00.30)
\drawline(54.30,00.30)(57.30,00.30)

\drawline(44.30,10.30)(47.30,10.30)
\drawline(49.30,10.30)(52.30,10.30)
\drawline(54.30,10.30)(57.30,10.30)

\drawline(44.30,20.30)(47.30,20.30)
\drawline(49.30,20.30)(52.30,20.30)
\drawline(54.30,20.30)(57.30,20.30)


\drawline(95.30,00.30)(98.30,00.30)
\drawline(100.30,00.30)(103.30,00.30)
\drawline(105.30,00.30)(108.30,00.30)

\drawline(95.30,10.30)(98.30,10.30)
\drawline(100.30,10.30)(103.30,10.30)
\drawline(105.30,10.30)(108.30,10.30)

\drawline(95.30,20.30)(98.30,20.30)
\drawline(100.30,20.30)(103.30,20.30)
\drawline(105.30,20.30)(108.30,20.30)


\put(76.00,10.00){\makebox(0,0)[cc]{$=$}}

\put(24.30,00.30){\makebox(0,0)[cc]{$p$}}
\put(24.30,10.30){\makebox(0,0)[cc]{$l$}}
\put(24.30,20.30){\makebox(0,0)[cc]{$k$}}
\end{picture}
\caption{An illustration of the equality
$l_{k,l}l_{p,l}=\varepsilon_k l_{p,l}$.}\label{fig:r5}
\end{figure}
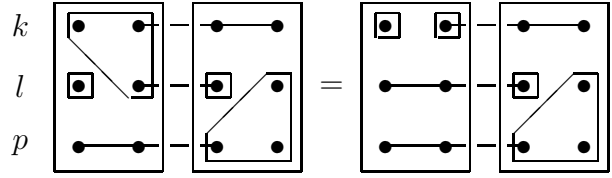

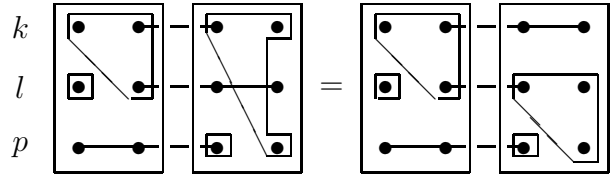
\begin{figure}
\special{em:linewidth 0.4pt} \unitlength 0.80mm
\linethickness{1pt}
\begin{picture}(150.00,35.00)

\put(34.00,00.00){\makebox(0,0)[cc]{$\bullet$}}
\put(34.00,10.00){\makebox(0,0)[cc]{$\bullet$}}
\put(34.00,20.00){\makebox(0,0)[cc]{$\bullet$}}
\put(44.00,00.00){\makebox(0,0)[cc]{$\bullet$}}
\put(44.00,10.00){\makebox(0,0)[cc]{$\bullet$}}
\put(44.00,20.00){\makebox(0,0)[cc]{$\bullet$}}

\drawline(30.00,-04.00)(30.00,24.00)
\drawline(30.00,24.00)(48.00,24.00)
\drawline(48.00,24.00)(48.00,-04.00)
\drawline(48.00,-04.00)(30.00,-04.00)

\drawline(34.30,00.30)(44.30,00.30)

\drawline(32.30,08.30)(32.30,12.30)
\drawline(32.30,12.30)(36.30,12.30)
\drawline(36.30,12.30)(36.30,08.30)
\drawline(36.30,08.30)(32.30,08.30)

\drawline(32.30,18.30)(32.30,22.30)
\drawline(32.30,22.30)(46.30,22.30)
\drawline(46.30,22.30)(46.30,08.30)
\drawline(46.30,08.30)(42.80,08.30)
\drawline(42.30,08.30)(32.30,18.30)


\put(57.00,00.00){\makebox(0,0)[cc]{$\bullet$}}
\put(57.00,10.00){\makebox(0,0)[cc]{$\bullet$}}
\put(57.00,20.00){\makebox(0,0)[cc]{$\bullet$}}
\put(67.00,00.00){\makebox(0,0)[cc]{$\bullet$}}
\put(67.00,10.00){\makebox(0,0)[cc]{$\bullet$}}
\put(67.00,20.00){\makebox(0,0)[cc]{$\bullet$}}

\drawline(53.00,-04.00)(53.00,24.00)
\drawline(53.00,24.00)(71.00,24.00)
\drawline(71.00,24.00)(71.00,-04.00)
\drawline(71.00,-04.00)(53.00,-04.00)

\drawline(55.30,-01.30)(55.30,02.30)
\drawline(55.30,02.30)(59.30,02.30)
\drawline(59.30,02.30)(59.30,-01.30)
\drawline(59.30,-01.30)(55.30,-01.30)

\drawline(57.30,10.30)(67.30,10.30)

\drawline(55.30,19.30)(55.30,22.30)
\drawline(55.30,22.30)(69.30,22.30)
\drawline(69.30,22.30)(69.30,18.30)
\drawline(69.30,18.30)(65.30,18.30)
\drawline(65.30,18.30)(65.30,02.30)
\drawline(65.30,02.30)(69.30,02.30)
\drawline(69.30,02.30)(69.30,-01.30)
\drawline(69.30,-01.30)(65.30,-01.30)
\drawline(65.30,-01.30)(55.30,19.30)


\put(85.00,00.00){\makebox(0,0)[cc]{$\bullet$}}
\put(85.00,10.00){\makebox(0,0)[cc]{$\bullet$}}
\put(85.00,20.00){\makebox(0,0)[cc]{$\bullet$}}
\put(95.00,00.00){\makebox(0,0)[cc]{$\bullet$}}
\put(95.00,10.00){\makebox(0,0)[cc]{$\bullet$}}
\put(95.00,20.00){\makebox(0,0)[cc]{$\bullet$}}

\drawline(81.00,-04.00)(81.00,24.00)
\drawline(81.00,24.00)(99.00,24.00)
\drawline(99.00,24.00)(99.00,-04.00)
\drawline(99.00,-04.00)(81.00,-04.00)

\drawline(85.30,00.30)(95.30,00.30)

\drawline(83.30,08.30)(83.30,12.30)
\drawline(83.30,12.30)(87.30,12.30)
\drawline(87.30,12.30)(87.30,08.30)
\drawline(87.30,08.30)(83.80,08.30)

\drawline(83.30,18.30)(83.30,22.30)
\drawline(83.30,22.30)(97.30,22.30)
\drawline(97.30,22.30)(97.30,08.30)
\drawline(97.30,08.30)(93.80,08.30)
\drawline(93.30,08.30)(83.30,18.30)


\put(108.00,00.00){\makebox(0,0)[cc]{$\bullet$}}
\put(108.00,10.00){\makebox(0,0)[cc]{$\bullet$}}
\put(108.00,20.00){\makebox(0,0)[cc]{$\bullet$}}
\put(118.00,00.00){\makebox(0,0)[cc]{$\bullet$}}
\put(118.00,10.00){\makebox(0,0)[cc]{$\bullet$}}
\put(118.00,20.00){\makebox(0,0)[cc]{$\bullet$}}

\drawline(104.00,-04.00)(104.00,24.00)
\drawline(104.00,24.00)(122.00,24.00)
\drawline(122.00,24.00)(122.00,-04.00)
\drawline(122.00,-04.00)(104.00,-04.00)

\drawline(106.30,08.30)(106.30,12.30)
\drawline(106.30,12.30)(120.30,12.30)
\drawline(120.30,12.30)(120.30,-02.30)
\drawline(120.30,-02.30)(116.30,-02.30)
\drawline(116.30,-02.30)(106.30,08.30)

\drawline(108.30,20.30)(118.30,20.30)

\drawline(106.30,-01.30)(106.30,02.30)
\drawline(106.30,02.30)(110.30,02.30)
\drawline(110.30,02.30)(110.30,-01.30)
\drawline(110.30,-01.30)(106.30,-01.30)


\drawline(44.30,00.30)(47.30,00.30)
\drawline(49.30,00.30)(52.30,00.30)
\drawline(54.30,00.30)(57.30,00.30)

\drawline(44.30,10.30)(47.30,10.30)
\drawline(49.30,10.30)(52.30,10.30)
\drawline(54.30,10.30)(57.30,10.30)

\drawline(44.30,20.30)(47.30,20.30)
\drawline(49.30,20.30)(52.30,20.30)
\drawline(54.30,20.30)(57.30,20.30)


\drawline(95.30,00.30)(98.30,00.30)
\drawline(100.30,00.30)(103.30,00.30)
\drawline(105.30,00.30)(108.30,00.30)

\drawline(95.30,10.30)(98.30,10.30)
\drawline(100.30,10.30)(103.30,10.30)
\drawline(105.30,10.30)(108.30,10.30)

\drawline(95.30,20.30)(98.30,20.30)
\drawline(100.30,20.30)(103.30,20.30)
\drawline(105.30,20.30)(108.30,20.30)


\put(76.00,10.00){\makebox(0,0)[cc]{$=$}}

\put(24.30,00.30){\makebox(0,0)[cc]{$p$}}
\put(24.30,10.30){\makebox(0,0)[cc]{$l$}}
\put(24.30,20.30){\makebox(0,0)[cc]{$k$}}
\end{picture}
\caption{An illustration of the equality
$l_{k,l}l_{k,p}=l_{k,l}l_{l,p}$.}\label{fig:r3}
\end{figure}

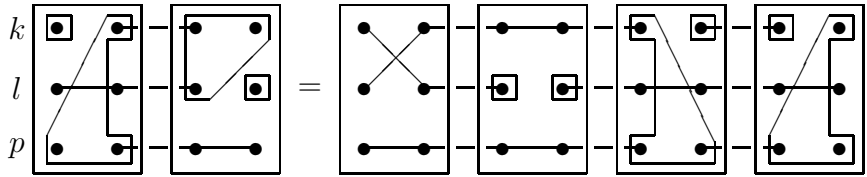
\begin{figure}
\special{em:linewidth 0.4pt} \unitlength 0.80mm
\linethickness{1pt}
\begin{picture}(150.00,35.00)
\put(11.00,00.00){\makebox(0,0)[cc]{$\bullet$}}
\put(11.00,10.00){\makebox(0,0)[cc]{$\bullet$}}
\put(11.00,20.00){\makebox(0,0)[cc]{$\bullet$}}
\put(21.00,00.00){\makebox(0,0)[cc]{$\bullet$}}
\put(21.00,10.00){\makebox(0,0)[cc]{$\bullet$}}
\put(21.00,20.00){\makebox(0,0)[cc]{$\bullet$}}

\drawline(07.00,-04.00)(07.00,24.00)
\drawline(07.00,24.00)(25.00,24.00)
\drawline(25.00,24.00)(25.00,-04.00)
\drawline(25.00,-04.00)(07.00,-04.00)

\drawline(11.30,10.30)(21.30,10.30)

\drawline(09.30,-02.30)(09.30,02.30)
\drawline(09.30,02.30)(19.30,22.30)
\drawline(19.30,22.30)(23.30,22.30)
\drawline(23.30,22.30)(23.30,18.30)
\drawline(23.30,18.30)(19.30,18.30)
\drawline(19.30,18.30)(19.30,02.30)
\drawline(19.30,02.30)(23.30,02.30)
\drawline(23.30,02.30)(23.30,-02.30)
\drawline(23.30,-02.30)(09.30,-02.30)

\drawline(09.30,18.30)(09.30,22.30)
\drawline(09.30,22.30)(13.30,22.30)
\drawline(13.30,22.30)(13.30,18.30)
\drawline(13.30,18.30)(09.30,18.30)


\put(34.00,00.00){\makebox(0,0)[cc]{$\bullet$}}
\put(34.00,10.00){\makebox(0,0)[cc]{$\bullet$}}
\put(34.00,20.00){\makebox(0,0)[cc]{$\bullet$}}
\put(44.00,00.00){\makebox(0,0)[cc]{$\bullet$}}
\put(44.00,10.00){\makebox(0,0)[cc]{$\bullet$}}
\put(44.00,20.00){\makebox(0,0)[cc]{$\bullet$}}

\drawline(30.00,-04.00)(30.00,24.00)
\drawline(30.00,24.00)(48.00,24.00)
\drawline(48.00,24.00)(48.00,-04.00)
\drawline(48.00,-04.00)(30.00,-04.00)

\drawline(34.30,00.30)(44.30,00.30)

\drawline(32.30,08.30)(32.30,22.30)
\drawline(32.30,22.30)(46.30,22.30)
\drawline(46.30,22.30)(46.30,18.30)
\drawline(46.30,18.30)(36.30,08.30)
\drawline(36.30,08.30)(32.30,08.30)

\drawline(42.30,08.30)(42.30,12.30)
\drawline(42.30,12.30)(46.30,12.30)
\drawline(46.30,12.30)(46.30,08.30)
\drawline(46.30,08.30)(42.30,08.30)


\put(62.00,00.00){\makebox(0,0)[cc]{$\bullet$}}
\put(62.00,10.00){\makebox(0,0)[cc]{$\bullet$}}
\put(62.00,20.00){\makebox(0,0)[cc]{$\bullet$}}
\put(72.00,00.00){\makebox(0,0)[cc]{$\bullet$}}
\put(72.00,10.00){\makebox(0,0)[cc]{$\bullet$}}
\put(72.00,20.00){\makebox(0,0)[cc]{$\bullet$}}

\drawline(58.00,-04.00)(58.00,24.00)
\drawline(58.00,24.00)(76.00,24.00)
\drawline(76.00,24.00)(76.00,-04.00)
\drawline(76.00,-04.00)(58.00,-04.00)

\drawline(30.00,-04.00)(30.00,24.00)
\drawline(30.00,24.00)(48.00,24.00)
\drawline(48.00,24.00)(48.00,-04.00)
\drawline(48.00,-04.00)(30.00,-04.00)

\drawline(62.30,00.30)(72.30,00.30)
\drawline(62.30,10.30)(72.30,20.30)
\drawline(62.30,20.30)(72.30,10.30)


\put(85.00,00.00){\makebox(0,0)[cc]{$\bullet$}}
\put(85.00,10.00){\makebox(0,0)[cc]{$\bullet$}}
\put(85.00,20.00){\makebox(0,0)[cc]{$\bullet$}}
\put(95.00,00.00){\makebox(0,0)[cc]{$\bullet$}}
\put(95.00,10.00){\makebox(0,0)[cc]{$\bullet$}}
\put(95.00,20.00){\makebox(0,0)[cc]{$\bullet$}}

\drawline(81.00,-04.00)(81.00,24.00)
\drawline(81.00,24.00)(99.00,24.00)
\drawline(99.00,24.00)(99.00,-04.00)
\drawline(99.00,-04.00)(81.00,-04.00)

\drawline(85.30,00.30)(95.30,00.30)
\drawline(85.30,20.30)(95.30,20.30)

\drawline(83.30,08.30)(83.30,12.30)
\drawline(83.30,12.30)(87.30,12.30)
\drawline(87.30,12.30)(87.30,08.30)
\drawline(87.30,08.30)(83.30,08.30)

\drawline(93.30,08.30)(93.30,12.30)
\drawline(93.30,12.30)(97.30,12.30)
\drawline(97.30,12.30)(97.30,08.30)
\drawline(97.30,08.30)(93.30,08.30)


\put(108.00,00.00){\makebox(0,0)[cc]{$\bullet$}}
\put(108.00,10.00){\makebox(0,0)[cc]{$\bullet$}}
\put(108.00,20.00){\makebox(0,0)[cc]{$\bullet$}}
\put(118.00,00.00){\makebox(0,0)[cc]{$\bullet$}}
\put(118.00,10.00){\makebox(0,0)[cc]{$\bullet$}}
\put(118.00,20.00){\makebox(0,0)[cc]{$\bullet$}}

\drawline(104.00,-04.00)(104.00,24.00)
\drawline(104.00,24.00)(122.00,24.00)
\drawline(122.00,24.00)(122.00,-04.00)
\drawline(122.00,-04.00)(104.00,-04.00)

\drawline(108.30,10.30)(118.30,10.30)

\drawline(116.30,18.30)(116.30,22.30)
\drawline(116.30,22.30)(120.30,22.30)
\drawline(120.30,22.30)(120.30,18.30)
\drawline(120.30,18.30)(116.30,18.30)

\drawline(106.30,-02.30)(106.30,02.30)
\drawline(106.30,02.30)(110.30,02.30)
\drawline(110.30,02.30)(110.30,18.30)
\drawline(110.30,18.30)(106.30,18.30)
\drawline(106.30,18.30)(106.30,22.30)
\drawline(106.30,22.30)(110.30,22.30)
\drawline(110.30,22.30)(120.30,01.30)
\drawline(120.30,01.30)(120.30,-02.30)
\drawline(120.30,-02.30)(106.30,-02.30)


\put(131.00,00.00){\makebox(0,0)[cc]{$\bullet$}}
\put(131.00,10.00){\makebox(0,0)[cc]{$\bullet$}}
\put(131.00,20.00){\makebox(0,0)[cc]{$\bullet$}}
\put(141.00,00.00){\makebox(0,0)[cc]{$\bullet$}}
\put(141.00,10.00){\makebox(0,0)[cc]{$\bullet$}}
\put(141.00,20.00){\makebox(0,0)[cc]{$\bullet$}}

\drawline(127.00,-04.00)(127.00,24.00)
\drawline(127.00,24.00)(145.00,24.00)
\drawline(145.00,24.00)(145.00,-04.00)
\drawline(145.00,-04.00)(127.00,-04.00)

\drawline(131.30,10.30)(141.30,10.30)

\drawline(129.30,-02.30)(129.30,02.30)
\drawline(129.30,02.30)(139.30,22.30)
\drawline(139.30,22.30)(143.30,22.30)
\drawline(143.30,22.30)(143.30,18.30)
\drawline(143.30,18.30)(139.30,18.30)
\drawline(139.30,18.30)(139.30,02.30)
\drawline(139.30,02.30)(143.30,02.30)
\drawline(143.30,02.30)(143.30,-02.30)
\drawline(143.30,-02.30)(129.30,-02.30)

\drawline(129.30,18.30)(129.30,22.30)
\drawline(129.30,22.30)(133.30,22.30)
\drawline(133.30,22.30)(133.30,18.30)
\drawline(133.30,18.30)(129.30,18.30)


\drawline(21.30,00.30)(24.30,00.30)
\drawline(26.30,00.30)(29.30,00.30)
\drawline(31.30,00.30)(34.30,00.30)

\drawline(21.30,10.30)(24.30,10.30)
\drawline(26.30,10.30)(29.30,10.30)
\drawline(31.30,10.30)(34.30,10.30)

\drawline(21.30,20.30)(24.30,20.30)
\drawline(26.30,20.30)(29.30,20.30)
\drawline(31.30,20.30)(34.30,20.30)


\drawline(72.30,00.30)(75.30,00.30)
\drawline(77.30,00.30)(80.30,00.30)
\drawline(82.30,00.30)(85.30,00.30)

\drawline(72.30,10.30)(75.30,10.30)
\drawline(77.30,10.30)(80.30,10.30)
\drawline(82.30,10.30)(85.30,10.30)

\drawline(72.30,20.30)(75.30,20.30)
\drawline(77.30,20.30)(80.30,20.30)
\drawline(82.30,20.30)(85.30,20.30)


\drawline(95.30,00.30)(98.30,00.30)
\drawline(100.30,00.30)(103.30,00.30)
\drawline(105.30,00.30)(108.30,00.30)

\drawline(95.30,10.30)(98.30,10.30)
\drawline(100.30,10.30)(103.30,10.30)
\drawline(105.30,10.30)(108.30,10.30)

\drawline(95.30,20.30)(98.30,20.30)
\drawline(100.30,20.30)(103.30,20.30)
\drawline(105.30,20.30)(108.30,20.30)


\drawline(118.30,00.30)(121.30,00.30)
\drawline(123.30,00.30)(126.30,00.30)
\drawline(128.30,00.30)(131.30,00.30)

\drawline(118.30,10.30)(121.30,10.30)
\drawline(123.30,10.30)(126.30,10.30)
\drawline(128.30,10.30)(131.30,10.30)

\drawline(118.30,20.30)(121.30,20.30)
\drawline(123.30,20.30)(126.30,20.30)
\drawline(128.30,20.30)(131.30,20.30)

\put(53.00,10.00){\makebox(0,0)[cc]{$=$}}

\put(04.30,20.30){\makebox(0,0)[cc]{$k$}}
\put(04.30,10.30){\makebox(0,0)[cc]{$l$}}
\put(04.30,00.30){\makebox(0,0)[cc]{$p$}}
\end{picture}
\caption{An illustration of the equality
$l_{p,k}r_{k,l}=s_{k,l}\varepsilon_l
r_{p,k}l_{p,k}$.}\label{fig:r8}
\end{figure}

Some examples of the relations satisfied by the generating
elements of the monoid $\PIP_n$, are given on Figures
\ref{fig:r5}, \ref{fig:r3} and \ref{fig:r8}.


Corollary~\ref{cor:can_form} and Proposition~\ref{prop:rel_pip_n}
imply that every element of $\PIP_n$ can be written as
$\varphi$-image of some canonical word from $S$. Now we are going
to show that such a presentation is unique.

\begin{theorem}
The map $\varphi:S\to\PIP_n$ from Proposition~\ref{prop:rel_pip_n}
is an isomorphism.
\end{theorem}

\begin{proof}
We are to prove that the map $\varphi$ is injective. Applying
Proposition~\ref{prop:canonical} it is enough to show that if
$\varphi$-images of values of two canonical words in $\PIP_n$ are
equal then these canonical words are equivalent. For this, we
compute the value of the image of a canonical word in $\PIP_n$.
For the word~\eqref{eq:canonical} this is the element
\begin{equation*}
\bigl\{\bigl(C_i\cup \sigma(D'_i)\bigr)_{1\leq i\leq s},
\{x\}_{x\in K_1}, \{\sigma(x')\}_{x\in K_2},
\{x,\sigma(x')\}_{x\in K_3}\bigr\},
\end{equation*}
where $K_1=M\cup (\cup_{i=1}^s(D_i\setminus C_i))$, $K_2=F=M\cup
(\cup_{i=1}^s(C_i\setminus D_i))$,
$K_3=X\setminus\bigl((\cup_{i=1}^s(C_i\cup D_i))\cup M\bigr)$. The
statement now follows from the definition of equivalent canonical
words.
\end{proof}

\vspace{0.3cm}

\noindent G.K.: Algebra, Department of Mathematics and Mechanics,
Kyiv Taras Shevchenko University, 64 Volodymyrska st., 01033 Kyiv,
UKRAINE, e-mail: {\tt akudr\symbol{64}univ.kiev.ua} \vspace{0.3cm}

\noindent V.M.: School of Mathematics and Statistics, University
of St Andrews, St Andrews, Fife, KY 16 9SS, SCOTLAND

\end{document}